\documentclass{amsart}
\usepackage{graphicx}
\graphicspath{{figures/}}
\usepackage{amssymb}
\usepackage[font={footnotesize,it}]{caption}
\usepackage{comment}
\usepackage{url}
\usepackage{hyperref}
\usepackage{caption}
\usepackage{mathtools}
\usepackage[colorinlistoftodos,prependcaption,textsize=tiny]{todonotes}

\numberwithin{theorem}{section}
\theoremstyle{remark}

\newcommand*\conj[1]{\overline{#1}}

\DeclarePairedDelimiter\floor{\lfloor}{\rfloor}
\DeclareMathOperator\asinh{asinh}

\begin{document}
\title[DE method for Riemann Zeta, Lerch and Dirichlet L-functions]{Double Exponential method for Riemann Zeta, Lerch and Dirichlet L-functions}
\author{Sandeep Tyagi}
\thanks{}
\address{377 Old Bedford Road, Luton, UK}
\email{tyagi\_sandeep@yahoo.com}

\begin{abstract}
We provide efficient methods to evaluate the Riemann zeta, the Lerch zeta and the Dirichlet $L$-functions. The method uses the Riemann-Siegel (RS) type formulas and a modified double exponential (MDE) quadrature method near the saddle point of appropriate integrands. We provide simplified derivations of the RS formulas, containing finite series sums and residual integrals, for the Lerch and the Dirichlet $L$-functions. The MDE method allows us to remove the contribution of the singularities near the contour of integration for the residual integrals. The method allows us to evaluate the residual integrals to any prescribed accuracy. The numerical cost of evaluating them is minimal compared to the main series sums. Thus even highly oscillatory integrands for these functions can be evaluated with the same complexity as the RS formula. In particular, the numerical complexity to evaluate  the $\zeta(s)$ and the Lerch zeta function with $s=\sigma+i t$ scales as  $\sqrt{t}$, and for Dirichlet $L$-function it scales as $\max(q, \sqrt{q t})$. The method allows the automatic setting of integral cutoffs and finite discretization size to achieve prescribed accuracy. Furthermore, it ensures that every halving of the discretization interval almost leads to doubling the accuracy of the results. Thus, it allows for a further check on the accuracy of the results obtained.

\end{abstract}

\maketitle

\section{Introduction}

Riemann zeta (RZ) function $\zeta(s)$ plays a vital role in analytical number theory and has applications in physics, probability theory, and statistics \cite{enwiki:1072466844, edwards1974riemann, titchmarsh1986theory}. Euler studied the function for real arguments, and Riemann extended it to cover the entire complex plane except $s=1$, where it has a simple pole. Riemann showed that the RZ function has a deep connection to the distribution of the prime numbers. Efficient numerical evaluation of the function is needed to check the validity of the Riemann hypothesis, obtaining RZ zeros to high accuracy \cite{platt2017isolating} and to check the mathematical connection between the moments of the zeta function to some models in random matrix theory \cite{hiary2011fast, rubinstein2005computational}.

The Lerch zeta function and the Dirichlet $L$-functions generalize the Riemann zeta function. Lerch function was introduced and investigated by Lerch \cite{lerch1887note}, and Lipschitz \cite{lipschitz1857untersuchung} in connection with Dirichlet's famous theorem on primes in arithmetic progression. The Lerch zeta function has been widely studied \cite{laurinvcikas2002lerch} and has found application in many branches of theoretical physics such as particle physics, thermodynamics, and statistical mechanics \cite{navasnumerical}. The Lerch zeta function covers many functions such as the RZ function, Hurwitz zeta function and Polylogarithms. The Dirichlet L-functions may be written as a linear combination of the Hurwitz zeta function at rational values. Thus, an efficient evaluation of the Lerch zeta function naturally leads to an efficient evaluation of the Dirichlet $L$-functions. However, due to immense interest in $L$-functions and the Generalized Riemann hypothesis, an efficient evaluation of the Dirichlet $L$-function is highly desired. Although the RS type formulas for Lerch zeta  \cite{balanzario2012riemann} and Dirichlet $L$-functions \cite{siegel1943contributions} do exist, their application to evaluating these functions to high accuracy has not been attempted. The error estimation with the RS formulas is very complex, even for the RZ function. Thus, in the absence of robust error bounds for RS formulas, one mainly depends on the Euler-Maclaurin method \cite{navasnumerical, johansson2015rigorous} for high accuracy computations. However, the Euler-Maclaurin formula is not suitable, especially if the imaginary part of the argument is large. In such cases, Euler-Maclaurin is extremely slow and puts a significant demand on computer memory requirements. For low accuracies, one may use the smoothed functional approach \cite{rubinstein2005computational, bailey2015crandall, dokchitser2004computing}. However, the smoothed functional technique is not suitable for high precision computations. It involves using incomplete gamma functions in its critical transition region, and there are no good ways to bound the error in this region.

This paper describes very accurate and highly efficient algorithms that apply the MDE method to the RS formulas for the RZ, the Lerch zeta function and the Dirichlet $L$-functions. The error estimation is straightforward, and one can achieve thousands of decimal point accuracy without the need to use exorbitant internal precision that is a hallmark of Euler-Maclaurin based approaches \cite{navasnumerical, johansson2015rigorous}. The method also does not require tedious error estimation based on heuristics \cite{de2011high, dereyna2022high}. For the Dirichlet $L$-functions, we give two implementations, namely that based on expanding the Dirichlet $L$-function in terms of Hurwitz zeta function and another one based on the Siegel formula \cite{siegel1943contributions}. The formulas are suitable for computing these functions for high values of $\Im(s)$, and the extra overhead required to calculate these functions to very high accuracy is typically negligible compared to the computation of principal series sums in the RS formulas.

The organization of the rest of the paper is as follows. In \autoref{sec:riemann}, \autoref{sec:lerch} and  \autoref{sec:dirichlet} we describe the three main functions investigated in this paper and derive the RS type formulas for two of them, namely the Lerch zeta function and the Dirichlet $L$-function. We then describe the MDE method in \autoref{sec:de}. In \autoref{sec:algo}  we develop algorithms for these functions. We report numerical tests in \autoref{sec:tests} and finally conclude with further scope of research in \autoref{sec:conclusion}. \\

{\bf Acknowledgements:} I would like to thank my wife Kirti and friend Adekunle Okusanya for constant encouragement and support to publish this work. I thank Prof. H. Cohen for pointing out a bug in Lerch function implementation and several other valuable suggestions.

\section{Riemann Zeta Function}
\label{sec:riemann}
The RZ function $\zeta (s)$ for $s=\sigma+i t$ is defined as $\zeta (s)=\sum_{n=1}^{\infty
}n^{-s}$. This sum converges for $\sigma >1$. The RZ function for all values of $\sigma$ can be defined by the analytic continuation of various integral forms representing the above sum. These integral representations remain valid for all $s$ except $s=1$, where the function has a simple pole. The RZ function is known to have non-trivial zeros only in the critical region $0 \le \sigma \le 1$. In particular, the Riemann hypothesis conjectures that all such zeros lie on the critical line $\sigma=1/2$. However, the conjecture remains unproven. 

There are numerical methods to compute this function, and ideally, one should be able to check this hypothesis numerically. For example, there is the Euler-Maclaurin formula, but the computational cost scales as $t$, and it requires very high internal precision for large $t$ values. Thus, the Euler-Maclaurin formula is only practical for $t$ values close to the real axis. On the other hand, the RS formula is the most useful of the currently known formulas for the RZ function for large $t$ computations. Although the computational cost of the RS formula scales as $\sqrt{t}$, it lacks rigorous and precise error bounds. 

Riemann derived the RS formula. The formula gives an integral representation of the function. One can make the contours of integrations pass through the integrands' saddle points, leading to principal series sums and residual integrals. The residual integrals are expressed in terms of a series expansion. While this expansion is good enough if low accuracy of the zeta function is required, its use for high accuracy computation is problematic because no error bounds have been established except for the particular case when $\sigma =1/2$ \cite{gabcke1979neue}. For the $\sigma=1/2$ case, Gabcke has derived a series of correction terms. Reyna \cite{de2011high} has generalized these results and has obtained expansions that allow one to evaluate the RZ function in the whole critical region.
However, some of the error estimations involve heuristics \cite{platt2017isolating}. Platt \cite{platt2017isolating} has used the bandwidth-limited method \cite{booker2006artin, molin2010integration, belabas2021numerical} to estimate the RZ function to high accuracy. However, these formulas lead to a computation cost proportional to $t$ for a single function evaluation at $s=\sigma+ i t$ as opposed to $\sqrt{t}$ complexity of the RS formula. Turing \cite{turing1953some} has given quadrature formulas for evaluating the RZ function, and Galway \cite{galway2001computing} improved these formulas. These methods do not provide good control over the error estimation. Lastly,  Rubinstein has given smoothed functional method to evaluate any general $L$-functions. It achieves $\sqrt{t}$ scaling. However, the explicit error bounds are difficult to obtain because there is no easy way to compute the incomplete gamma function to a prescribed accuracy, especially in its transition region.

We discuss two of the above methods in more detail. First we consider the smoothed functional method to compute the RZ function. The method has been considerd by a number of authors such as Crandall \cite{bailey2015crandall} and Rubinstein \cite{rubinstein2005computational}. It has a computational complexity of $t^{1/2}$ \cite{rubinstein2005computational}. There are similar methods for the computation of the Lerch zeta \cite{balanzario2012riemann} and Dirichlet $L$-function \cite{siegel1943contributions}  as well. The formula is given by
\begin{equation}  \label{smoothed approx}
\pi^{-s/2}\Gamma(s/2)\zeta(s)\delta^{-s} = -\frac{1}{s}-\frac{\delta^{-1}}{%
1-s} +\sum_{n=1}^{\infty} G(s/2,\pi n^2\delta^2)
+\delta^{-1}\sum_{n=1}^{\infty} G((1-s)/2,\pi n^2/\delta^2),
\end{equation}
where $G(z,w)$ is defined in terms of the incomplete Gamma function $\Gamma(z,w)$ as 
$G(z,w) =w^{-z}\Gamma(z,w) = \int_1^{\infty} e^{-wx}x^{z-1}\,dx$, $\Re(w) > 0$, and 
$\delta$ is a free complex parameter of modulus one. A value of $\delta=1$ leads to a situation where the evaluation of $\zeta(\sigma+it)$ for very large $t$ leads to a catastrophic cancellation on the right hand side.  As a result, it requires working with very high internal precision which is not practical for large $t$ values. On the other hand if the $\delta$ is chosen so as to counterbalance this decay of the Gamma function, one can avoid the cancellation on the RHS. However, this requires calculation of the incomplete Gamma function for complex $z$ and $w$. The calculation of this function in the transition region where $z \approx w$ is a challenging problem \cite{rubinstein2005computational}.  

The second method is based on RS integral formula \cite{galway2001computing, edwards1974riemann}
\begin{equation}
\zeta \left( s\right) =I_{0}\left( s\right) +\chi \left( s\right) \overline{%
I_{0}\left( 1-\overline{s}\right) },  \label{z2}
\end{equation}%
where the bar denotes the complex conjugation, 
\begin{equation}
\chi(s) = \pi^{s-1/2} \frac{\Gamma\left(\frac{1-s}{2}\right)}{\Gamma\left(\frac{s}{2}\right)},
\end{equation}
\begin{eqnarray}
I_{0}\left( s\right)  &=&\int_{0\swarrow 1}g\left( z\right) dz
\label{I0},
\end{eqnarray}
and
\begin{eqnarray}
g(z)&=&\frac{e^{i\pi z^{2}}}{e^{i\pi z} -e^{-i\pi z}} z^{-s}.
\end{eqnarray}%
The integration contour is chosen as a straight line from the first quarter to the third quarter. The contour intersects the x-axis between $0$ and $1$. The contour in $I_0$ can be moved on the x-axis to pass near the saddle point of the integrand. This displacement leads to a series sum and a residual integral. The residual integral can be expressed in the form of a series \cite{hiary2016alternative}, and overall we obtain the following expansion for the RZ function on the critical line:
\begin{equation}  
\label{eq:rsform}
e^{i\theta(t)}\zeta(1/2+it) = 2\, \Re \left(e^{-i\theta(t)}\sum_{n=1}^{n_1} 
\frac{e^{it \log n}}{\sqrt{n}}\right) - \frac{(-1)^{n_1}}{\sqrt{a}}
\sum_{r=0}^m \frac{C_r(z)}{a^r} + R_m(t),
\end{equation} 
where
\begin{equation}
e^{i \theta(t)} =\left( \frac{\Gamma\left(\frac{1}{4}+i t\right)}{\Gamma\left(\frac{1}{4}-i t\right)}\right)^{\frac{1}{2}} \pi^{-i \frac{t}{2}},
\end{equation}
$t>2\pi$, $a=\sqrt{t/(2\pi)}$, $n_1=\lfloor
a\rfloor $ is the integer part of $a$ and $z= 1-2(a-\lfloor a\rfloor)$. 
The $C_r(z)$ are given in terms of the derivatives of $F(z)=\cos ((\pi/2)(z^2+3/4))(\cos(\pi z))^{-1}$.
Gabcke \cite{gabcke1979neue}  has derived $C_r(z)$ on the critical line $\sigma=1/2$ and has obtained explicit error bounds for $R_m(t)$ assuming $t>200$. This work has been extended by Reyna for $\sigma \in (0,1)$. The  RZ implementation in MPMATH package \cite{johansson2010mpmath} is based on this implementation. However, some of the calculations depend on heuristics \cite{de2011high, de2013exact}. Note that both Gabcke and Reyna's formulas use the series expansion above. However, Berry has shown that one cannot obtain better and better accuracy by just continuing this series as the series is divergent. The series needs to be stopped at an optimal number of terms.

Starting with $I_0$ in \autoref{eq:I0}, Galway \cite{galway1999fast} has developed quadrature formulas for RZ function, but they lack in sharp and rigorous error bounds. Finally, there are methods pioneered by Hairy \cite{hiary2011fast} that try to compute the series sum over $N=\floor{\sqrt{t/(2 \pi)}}$ terms faster by splitting it into smaller sums which are then transformed into truncated theta sums \cite{hiary2011nearly}. The best scaling that one can achieve this way is $t^{4/13}$. However, except for Reyna's method, none of the methods above can calculate the RZ function to arbitrary accuracy in a reasonable amount of time.

Thus, overall it appears that there is no way to compute the RZ function in the critical region rigorously, especially for large $t$ that can also scale as $\sqrt{t}$ or better. We address this problem in this paper. First, we derive formulas for the RZ function in the critical region by applying the MDE method to the residual integrals of the RS formula. Our method leads to formulas that show $\sqrt{t}$ behaviour and, in addition, provide control over the error. Furthermore, we show that the accuracy obtained from the formula doubles with every halving of the interval $h$ used to approximate the integrals. In addition, we validate the work of Reyna by showing that we did not encounter even a single case where our results differed from his even to as high accuracy as $1000$ decimal places.

\section{Lerch Zeta function}
\label{sec:lerch}
\begin{figure}[ht]
\centering
\includegraphics[width=0.5\textwidth]{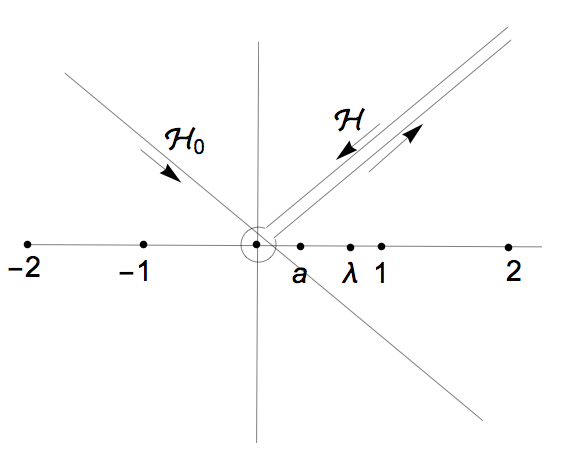}
\caption{The figure shows the Hankel contour $\mathcal{H}$. The branch cut is at $\pi/4$. Thus angle just past this branch cut should be taken to be $-3\pi/4$.}
\label{fig:f1}
\end{figure}

Lerch zeta function is defined as:
\begin{align}
L(s,\lambda,a)&=\sum_{n=0}^{\infty} \frac{e^{2i \pi \lambda n}}{(n+a)^s},
\label{lerch1a} 
\end{align}
where $s=\sigma + i t$,  $\sigma>1$ and we consider the case where both  $a$ and $\lambda$ are real with values between 0 and 1. We believe the extension of the formulas to complex $a$ and $\lambda$ should be straightforward.  Following Balanzario \cite{balanzario2012riemann},  we can obtain a representation of the Lerch function valid in the entire complex plane for $\sigma$. We provide an easy derivation of this representation. To start with, we first multiply both sides of \autoref{lerch1a} by $\Gamma(s)$. This leads to
\begin{align}
L(s,\lambda,a)&= \frac{1}{a^s}+\frac{(2\pi)^s}{\Gamma(s)}\int_{0}^{\infty}\frac{z^{s-1} e^{-2\pi a z}}{e^{2\pi z -2\pi \lambda i }-1} dz \nonumber \\
&= \frac{1}{a^s}+\frac{(2\pi)^s}{\Gamma(s)}\int_{0}^{\infty e^{-i \pi/4}}\frac{z^{s-1} e^{-2\pi a z}}{e^{2\pi z -2\pi \lambda i }-1} dz,
\label{lerch1}
\end{align}
where moving the integration contour from $(0,\infty)$ to $(0, \infty e^{-i \pi/4})$ is justified as the integral contribution from $(\infty, \infty e^{-i\pi/4})$ is zero and no residues are encountered when moving the contour.  We start with the identity,
\begin{align}
\int_{\mathcal{L}} \frac{e^{i \pi w^2+2 i \pi z w+\pi w}}{e^{2 \pi w}-1} dw & = \frac{i}{e^{2 \pi z}-1}-\frac{i e^{-i \pi z^2+\pi z}}{e^{2 \pi z}-1},
\label{eq:id1}
\end{align}
where $\mathcal{L}$ represents a straight line from the third quarter to the first and intersecting the real $x$-axis between 0 and 1.  The $\mathcal{L}$ can be parameterised as $-\alpha i+\xi e^{i \pi/4}$ with $\alpha \in (0,1)$ and $\xi$ a real number. Rearranging the terms in \autoref{eq:id1} and substituting $z$ with $(z-i \lambda)$, and $w$ with $i(w-a)$ yields
\begin{align}
\frac{1}{e^{2 \pi (z- i \lambda)}-1} = -i \int_{\mathcal{L}/i+a} \frac{e^{-i \pi (w-a)^2- 2\pi (z-i \lambda) (w-a)+\pi (w-a) }}{e^{2 i \pi (w-a)}-1} dw +\frac{e^{-i \pi (z- i \lambda)^2 +\pi (z-i \lambda)}}{e^{2 \pi( z-i \lambda)}-1}.
\label{lerch3}
\end{align}
Starting with \autoref{lerch3}, we multiply both sides of the equation by $z^{s-1} e^{-2 \pi a z}$ and integrate over $z$ from $0$ to $\infty \exp(-i \pi/4)$. This, using \autoref{lerch1},  leads to
\begin{align}
\frac{\Gamma(s) }{(2\pi)^s}(L(s,\lambda,a)-a^{-s})&= - i \int_{\mathcal{L}/i+a} \frac{e^{-i \pi (w- a)^2+ 2 i \pi \lambda (w-a)+i \pi (w-a) }}{e^{2 i \pi (w-a)}-1} dw \int_{0}^{\infty e^{-i \pi/4}} e^{2 i \pi z w} z^{s-1} dz \\
&+\int_{0}^{\infty e^{-i \pi/4}} \frac{e^{-i \pi (z-i \lambda)^2 + \pi (z-i \lambda)-2 \pi a z}}{e^{2 \pi (z- i \lambda)}-1} z^{s-1} dz.
\end{align}
Using the definition of the Gamma function,
\begin{align}
\int_{0}^{\infty e^{-i \pi/4}} e^{2 i \pi z w } z^{s-1} dz  = e^{i \pi s/2} (2\pi w)^{-s} \Gamma(s),
\label{eq:gamma}
\end{align}
we obtain
\begin{align*}
\frac{\Gamma(s) }{(2\pi)^s}(L(s,\lambda,a)-a^{-s})&= - i e^{i \pi s/2} (2\pi)^{-s} \Gamma(s) \int_{\mathcal{L}/i+ a} \frac{e^{-i \pi (w-a)^2+ 2 i \pi \lambda (w- a)+i \pi (w- a) }}{e^{2 i \pi (w-a)}-1} w^{-s} dw \\
&+ \int_{0}^{\infty e^{-i \pi/4}}\frac{e^{-i \pi (z-\lambda i)^2 + \pi (z-\lambda i)-2 \pi a z}}{e^{2 \pi z-2 \pi \lambda i}-1} z^{s-1} dz.
\end{align*}
In the second integral we replace  $z$ with $-w i$ and obtain
\begin{align*}
\frac{\Gamma(s) }{(2\pi)^s}(L(s,\lambda,a)-a^{-s})&= e^{i \pi s/2} (2\pi)^{-s} (i)^{-s} \Gamma(s) \int_{\mathcal{L}/i+a} \frac{e^{-i \pi (w-a)^2+ 2 i \pi \lambda (w-a)+i \pi (w- a) }}{e^{2 \pi  i (w-a)}-1} w^{-s} dw \\
&+  (-i)^{s} \int_{0}^{\infty e^{-i \pi/4}/(-i)}\frac{e^{i \pi (w+\lambda)^2 -i \pi (w+\lambda)+2 i \pi a w}}{e^{-2 i \pi (w+\lambda)}-1} w^{s-1} dw.
\end{align*}
This leads to 
\begin{align}
L(s,\lambda,a)&= a^{-s} + e^{-i \pi a (1+a+2\lambda)}\int_{\mathcal{L}/i+a} \frac{e^{-i \pi w^2+ 2i \pi w (a+\lambda)}}{e^{\pi  i w- 2 i \pi a}-e^{-i \pi w}} w^{-s} dw \nonumber \\
&- \frac{e^{-i\pi s/2}}{(2\pi)^{-s} \Gamma(s)} e^{i \pi \lambda (\lambda+1)}\int_{0}^{\infty e^{i \pi/4}}\frac{e^{i \pi w^2 +2 i \pi w(a+\lambda)}}{e^{i \pi w+2 \pi \lambda i}-e^{-i \pi w}} w^{s-1} dw.
\label{eq:l2}
\end{align}
As the contour $\mathcal{L}$ is parameterised as $- i \alpha+\xi e^{i \pi/4}$, the contour $\mathcal{L}/i+a$ corresponds to $(a-\alpha) + \xi e^{-i \pi/4}$ which represents a contour from the second quarter to the fourth quarter. The additive part here is in the range $(a-1, a)$, but to avoid the branch cut at the negative $x$-axis, we restrict the range to region $(0, a)$. Noting that the first integral on the right-hand side of \autoref{eq:l2} has a pole at $w=a$, we can displace the contour parallel to itself along the positive $x$-axis only up to $a$ without encountering a pole. The contour can be displaced further to the right after taking the contribution of the poles crossed into account.
Denoting
\begin{align*}
f(w)&= \frac{e^{-i \pi w^2+ 2i \pi w (a+\lambda)}}{e^{\pi  i w- 2 i \pi a}-e^{-i \pi w}} w^{-s} dw, \\
\end{align*}
the first integral on the RHS of \autoref{eq:l2} can be written as
\begin{align*}
I_1 &=e^{-i \pi a (1+a+2\lambda)}\int_{0 \searrow a} f(w) dw
\end{align*}
and displacing the contour to the right over the pole at $a$, we have
\begin{align*}
I_1 &=a^{-s} + \int_{a \searrow (a+1)} f(w) dw.
\end{align*}
To simplify the second part on the RHS of \autoref{eq:l2} we denote 
\begin{align*}
g(w)&= \frac{e^{i \pi w^2 +2 i \pi w(a+\lambda)}}{e^{i \pi w+2 \pi \lambda i}-e^{-i \pi w}},
\end{align*}
and note
\begin{align*}
\int_{0}^{\infty e^{i \pi/4}} g(w) w^{s-1} dz &= \frac{1}{1-e^{-2i \pi (s-1)}} \int_{\mathcal{H}} g(w) w^{s-1} dw \\
&=\frac{e^{i \pi s}}{2 i \sin(\pi s)} \int_{\mathcal{H}} g(w) w^{s-1} dw,
\end{align*}
where contour $\mathcal{H}$ is shown in \autoref{fig:f1}. In particular, this contour has a negative sign off compared to the same contour defined in Balanzario \cite{balanzario2012riemann}.
Using the relations above, we can express $L(s,\lambda,a)$ in its final form as
\begin{align}
L(s,\lambda,a)&=  e^{-i \pi a (1+a+2\lambda)}\int_{0 \searrow a} \frac{e^{-i \pi w^2+ 2i \pi w (a+\lambda)}}{e^{\pi  i w- 2 i \pi a}-e^{-i \pi w}} w^{-s} dw \\ \nonumber
&+ i e^{i \pi \lambda (\lambda+1)+i \pi s/2} \frac{\Gamma(1-s)}{(2\pi)^{1-s}}\int_{\mathcal{H}}\frac{e^{i \pi w^2 +2 i \pi w(a+\lambda)}}{e^{i \pi w+2 \pi \lambda i}-e^{-i \pi w}} w^{s-1} dw.
\label{eq:lerch}
\end{align}
By analytic continuation, \autoref{eq:lerch} remains valid in the entire complex plane $s$.

\subsection{Hurwitz and Riemann Zeta Functions:}
\label{subsec:hurwitz}
From the above representation of the Lerch zeta function, we can obtain an expression for the Hurwitz zeta function by taking the limit $\lambda \to 0$:
\begin{align*}
\zeta(s, a)&=  e^{-i \pi a (1+a)}\int_{0 \searrow a} \frac{e^{-i \pi w^2+ 2i \pi w a}}{e^{\pi  i w- 2 i \pi a}-e^{-i \pi w}} w^{-s} dw \\
&+i e^{i \pi s/2} \frac{\Gamma(1-s)}{(2\pi)^{1-s}}\int_{\mathcal{H}}\frac{e^{i \pi w^2 +2 i \pi w a}}{e^{i \pi w}-e^{-i \pi w}} w^{s-1} dw.
\end{align*}
We can also obtain a representation for $\zeta(s)$ by taking the limit $a \to 0$:
\begin{align*}
\zeta(s)&=\lim _{a \to 0} \sum_{n=0}^{\infty} \frac{1}{(n+a)^s}-\frac{1}{a^s}, \\
             &=\int_{0 \searrow 1 } \frac{e^{-i \pi w^2}}{e^{i \pi w}-e^{-i \pi w}} w^{-s} dw
            + i e^{\frac{\pi s}{2}i } \frac{\Gamma(1-s)}{(2\pi)^{1-s}} \int_{\mathcal{H}} \frac{e^{i \pi w^2}}{e^{i \pi w}-e^{-i \pi w}} w^{s-1} dw.
\end{align*}
Using the relation
\begin{align}
\pi^{-\frac{s}{2}} \zeta(s) \Gamma\left(\frac{s}{2}\right) = \pi^{-\frac{1-s}{2}} \zeta(1-s) \Gamma\left(\frac{1-s}{2}\right),
\end{align}
and the fact that the integrals around the Hankel contour can be simplified for the case of $a=0$ and $\lambda=0$, we have
\begin{align}
\int_{\mathcal{H}} g(w) w^{s-1} dw &= \frac{1-e^{-2i \pi (s-1)}}{1+e^{-i \pi s}} \int_{0 \nearrow 1} g(w) w^{s-1} dw \nonumber \\
&=-(1-e^{-i \pi s})\int_{0 \swarrow 1} g(w) w^{s-1} dw,
\label{eq:change}
\end{align}
where we have used $g(-w)=e^{i \pi} g(w)$ which is true on account of $a$ and $\lambda$ both being zero. Note that we have changed the direction of the integration contour in \autoref{eq:change}. With these changes we obtain
\begin{align}
\zeta(s)  &= \pi^{s-\frac{1}{2}}  \frac{\Gamma\left(\frac{1-s}{2}\right)}{\Gamma\left(\frac{s}{2}\right)} \zeta(1-s) \nonumber \\
&=\pi^{s-\frac{1}{2}}  \frac{\Gamma\left(\frac{1-s}{2}\right)}{\Gamma\left(\frac{s}{2}\right)} \left( \int_{0 \searrow 1 } \frac{e^{-i \pi w^2}}{e^{i \pi w}-e^{-i \pi w}} w^{s-1} dw
            -e^{\frac{-\pi s}{2}i } \frac{\Gamma(s)}{(2\pi)^{s}} \int_{\mathcal{H}} \frac{e^{i \pi w^2}}{e^{i \pi w}-e^{-i \pi w}} w^{-s} dw\right) \nonumber \\
&=\pi^{s-\frac{1}{2}} \frac{\Gamma\left(\frac{1-s}{2}\right)}{\Gamma\left(\frac{s}{2}\right)}  \int_{0 \searrow 1 } \frac{e^{-i \pi w^2}}{e^{i \pi w}-e^{-i \pi w}} w^{s-1} dw + q(s) \int_{0 \swarrow 1} \frac{e^{i \pi w^2}}{e^{i \pi w}-e^{-i \pi w}} w^{-s} dw \nonumber,
\end{align}
where
\begin{align}
q(s)=\left(e^{\frac{-\pi s}{2}i } + e^{\frac{\pi s}{2}i }\right) \pi^{s-\frac{1}{2}} \frac{\Gamma\left(\frac{1-s}{2}\right)}{\Gamma\left(\frac{s}{2}\right)}   \frac{\Gamma(s)}{(2\pi)^{s}}.
\end{align}
Now using the relations
\begin{align}
\Gamma(z) \Gamma(1-z) &= \frac{\pi}{\sin(\pi s)},
\end{align}
and
\begin{align}
\Gamma(z) \Gamma(z+\frac{1}{2})=2^{1-2z} \sqrt{\pi} \Gamma(2z),
\end{align}
it is easy to show that  $q(s)=1$ and thus this leads to well known RZ expansion as a special case of the Lerch zeta function:
\begin{align}
\zeta(s)  &= \int_{0 \swarrow 1} \frac{e^{i \pi w^2}}{e^{i \pi w}-e^{-i \pi w}} w^{-s} dw+\pi^{s-\frac{1}{2}} \frac{\Gamma\left(\frac{1-s}{2}\right)}{\Gamma\left(\frac{s}{2}\right)}  \int_{0 \searrow 1 } \frac{e^{-i \pi w^2}}{e^{i \pi w}-e^{-i \pi w}} w^{s-1} dw.
\end{align}

\section{Dirichlet $L$-Function: Siegel Method}
\label{sec:dirichlet}
The Dirichlet $L$-function is expressible in terms of the Hurwitz function as
\begin{align}
L(s,\chi) &= \sum_{n=1}^{\infty}\frac{\chi(n)}{n^s} =\frac{1}{m^s}\sum_{n=1}^{m}{\chi(n)}\zeta\left(s, \frac{n}{m}\right),
\label{eq:dirichlet_hurwitz}
\end{align}
where $\chi(n)$ denotes the Dirichlet character with period $m$.
With the formula above, the  cost of computing the function scales as $q\sqrt{t/(2\pi)}$. However, 
Siegel \cite{siegel1943contributions} has derived a formula for which the computational complexity scales as $\sqrt{q t/(2\pi)}$. We provide an easy proof of this formula and then show how one can use it to evaluate the Dirichlet $L$- function to any prescribed accuracy with scaling as $\min(\sqrt{q t/(2\pi)},q)$.  The formula is given as
\begin{equation}
\rho \left(\frac{m}{\pi} \right)^{s/2} \Gamma\left(\frac{s+a}{2}\right) L(s,\chi) = \mu(s) + \overline{\mu(1-\overline{s})},
\label{eq:dirichlet}
\end{equation}
where 
\begin{equation}
\mu(s) = \rho \left(\frac{m}{\pi} \right)^{s/2} \Gamma\left(\frac{s+a}{2}\right)  \lambda(s),
\label{eq:mu}
\end{equation}
\begin{align}
\lambda(x)= \frac{1}{2 \pi i} \int_{0 \swarrow 1} e^{\frac{i \pi}{m} x^2} x^{-s} W(x) dx,
\label{eq:lambda}
\end{align}
\begin{align}
W(x)=\frac{\pi}{2m} \sum_{n=1}^{2m} \chi(n) e^{-i \pi \frac{n^2}{m}} \cot\left( \pi \frac{x-n}{2m}\right),
\label{eq:w}
\end{align}
\begin{equation}
\rho = i^{-a/2} C^{-1/2} m^{1/4},
\end{equation}
\begin{equation}
C = \sum_{n=1}^{m}  \chi(n) e^{-2 \pi i \frac{n}{m}} ,
\end{equation}
and 
\begin{equation}
a = \frac{1-\chi(-1)}{2}.
\label{eq:a}
\end{equation}
The main element to prove the Siegel formula in \autoref{eq:mu} is the identity
\begin{align}
\int_{0 \swarrow 1} e^{\frac{i \pi}{m} x^2 - \frac{2 i \pi}{m} \xi x} W(x) dx  &= \frac{2 \pi i}{1-e^{-2i \pi \xi}}\sum_{n=1}^{m} \chi(n) e^{\frac{2i \pi n}{m} \xi} - C e^{-\frac{i \pi}{m}\xi^2} \conj{W(\conj{\xi})}.
\label{e1}
\end{align}
We provide proof of the identity in \autoref{e1} that is different to the one provided by Siegel. We start by writing $W(x)$ in \autoref{eq:w} as
\begin{align}
W(x)=\frac{\pi i}{m} \sum_{n=1}^{2m} \chi(n) \frac{ e^{-i \pi \frac{n^2}{m}}}{1-e^{-i \pi \frac{x-n}{m}}} - \frac{\pi i}{2m} \sum_{n=1}^{2m} \chi(n)  e^{-i \pi \frac{n^2}{m}}.
\end{align}
To prove it we start with integral
\begin{align}
J_{0 \swarrow 1} = \int_{0 \swarrow 1} \frac{e^{q i \pi (x-\xi)^2} }{1-e^{-2i \pi x}} dx,
\label{eq:j01}
\end{align}
and applying residue theorem note that
\begin{align}
J_{-1 \swarrow 0} &= \int_{-1 \swarrow 0} \frac{e^{q i \pi (x-\xi)^2} }{1-e^{-2i \pi x}} dx\\ \label{eq:jint}
&= J_{0 \swarrow 1} + e^{q i \pi \xi^2}.
\end{align}
On the other hand, substituting $x+1$ in place of $x$ in \autoref{eq:j01}, we have
\begin{align}
J_{0 \swarrow 1} &= \int_{-1 \swarrow 0} \frac{e^{q i \pi (x+1-\xi)^2}}{1-e^{-2i \pi (x+1)}} dx \nonumber \\
&= e^{q i \pi (1-2\xi)} \int_{-1 \swarrow 0} e^{q i \pi (x-\xi)^2} \frac{e^{2 q i \pi x}}{1-e^{-2 i \pi x}} dx \nonumber \\
&=  e^{q i \pi (1-2\xi)} \left(J_{-1 \swarrow 0}  + \int_{-1 \swarrow 0} e^{q i \pi (x-\xi)^2}  \frac{e^{2 q i \pi x}-1}{1-e^{-2i \pi x}} dx \right).
\label{eq:j01a}
\end{align}
The integral on the RHS of \autoref{eq:j01a} can be calculated using the relations
\begin{align}
\frac{e^{2 q i \pi x}-1}{1-e^{-2i \pi x}} &=\sum_{n=1}^{q} e^{2 i \pi n x},
\end{align}
and
\begin{align}
\int_{-1 \swarrow 0} e^{q i \pi (x-\xi)^2}  e^{2 q i \pi k x} dx &= e^{-i \pi \frac{k^2}{q}+2 \pi k i \xi} \int_{-1 \swarrow 0} e^{q i \pi (x-\xi+\frac{k}{q})^2}  dx \nonumber \\
&=-e^{i \pi/4} \sqrt{\frac{1}{q}} e^{-i \pi \frac{k^2}{q}+2 \pi k i \xi}.
\end{align}
We thus obtain
\begin{align}
J_{-1 \swarrow 0} (1-e^{-q i \pi (1-2\xi)}) &=   -e^{-q i \pi (1-2\xi)} e^{q i \pi \xi^2}+e^{i \pi/4} \sqrt{\frac{1}{q}} \sum_{k=1}^q  e^{-i \pi \frac{k^2}{q}+2 \pi k i \xi}.
\end{align}
We now replace $q$ with $4m$, $x$ with $(x-n)/(2m)$ and $\xi$ with $(\xi-n)/(2m)$ in \autoref{eq:jint} and denote by $J(\xi,n)$ the integral
\begin{align}
J(\xi,n)&= \frac{1}{2m} \int_{-1 \swarrow 0} \frac{e^{ \frac{i \pi}{m} (x-\xi)^2}}{1-e^{-i \pi \frac{(x-n)}{m}}} dx.
\end{align}
This gives
\begin{align}
({1-e^{4 i \pi \xi}}) J(\xi,n)  &= -e^{4 i \pi \xi} e^{ \frac{i \pi}{m} (\xi-n)^2} +e^{i \pi/4} \sqrt{\frac{1}{4m}} \sum_{k=1}^{4m}  e^{-i \pi \frac{k^2}{4m}+\pi k i \frac{\xi-n}{m}}.
\label{eq:jpsi}
\end{align}
We multiply both sides of \autoref{eq:jpsi} with $\chi(n)e^{-i \pi \frac{n^2}{m}}$ and sum over $n$ from $1$ to $2m$ to obtain
\begin{align}
({1-e^{4 i \pi \xi}}) \sum_{n=1}^{2m} J(\xi,n) \chi(n) e^{-i \pi \frac{n^2}{m}} &= J_1+J_2,
\end{align}
where
\begin{align}
J_1 &= -e^{4 i \pi \xi}\sum_{n=1}^{2m} \chi(n) e^{-i \pi \frac{n^2}{m}} e^{ \frac{i \pi}{m} (\xi-n)^2} \nonumber \\
&=-e^{4 i \pi \xi} e^{ \frac{i \pi}{m} \xi^2} (1+e^{-2i \pi \xi}) \sum_{n=1}^{m} \chi(n) e^{-2i \pi \frac{\xi n}{m}},
\end{align}
and
\begin{align}
J_2=e^{i \pi/4} \sqrt{\frac{1}{4m}} \sum_{k=1}^{4m} e^{-i \pi \frac{k^2}{4m}+\pi k i \frac{\xi}{m}} \sum_{n=1}^{2m} \chi(n) e^{-i \pi \frac{n^2}{m}-\pi k i \frac{n}{m}}.
\end{align}
To simplify $J_2$, we use the identity
\begin{align}
\chi(n) = \frac{1}{\conj{C} } \sum_{k'=1}^{m} \conj{\chi}(k') e^{2i \pi \frac{k' n}{m}},
\end{align}
to derive
\begin{align}
\sum_{n=1}^{2m}\chi(n)e^{-i \pi \frac{n^2}{m}- i \pi \frac{k n}{m}} &=\frac{1}{ \conj{C}} \sum_{k'=1}^{m} \conj{\chi}(k') \sum_{n=1}^{2m} e^{-i \pi \frac{n^2}{m} - i \pi \frac{(k-2k') n}{m}} \nonumber \\
&= \frac{1}{ \conj{C}} \sum_{k'=1}^{m} \conj{\chi}(k') \sum_{n=1}^{m} e^{-i \pi \frac{n^2}{m} - i \pi \frac{(k-2k') n}{m}}  (1+e^{-i \pi (m+k)})
\end{align}
Now we can use the fact that if $(m+k)$ is even then
\begin{align}
\sum_{n=1}^{m}e^{-i \pi \frac{n^2}{m}-i \pi \frac{k n}{m}} &= e^{\frac{i \pi}{4}\left(\frac{k^2}{m}-1\right)}\sqrt{m}.
\end{align}
Thus, we obtain after substituting $k$ with $(k-2k\rq{})$
\begin{align}
 \sum_{n=1}^{2m}\chi(n)e^{-i \pi \frac{n^2}{m}- i \pi \frac{k n}{m}} 
&=\frac{\sqrt{m}}{ \conj{C}} e^{-\frac{i \pi}{4}} \sum_{k'=1}^{m} \conj{\chi}(k') e^{\frac{i \pi}{4}\left(\frac{(k-2k')^2}{m}\right)} (1+e^{-i \pi (m+k)}) \nonumber \\
&=\frac{\sqrt{m}}{ \conj{C}} e^{-\frac{i \pi}{4}} e^{i \pi \frac{k^2}{4m}} \sum_{k'=1}^{m} \conj{\chi}(k') e^{i \pi \frac{k'^2}{m} }e^{-\frac{i \pi k k'}{m}}(1+e^{-i \pi (m+k)}).
\end{align}
This gives
\begin{align}
J_2&=\frac{1-e^{4 i \pi \xi}}{2 \conj{C}} \sum_{k'=1}^{m} \conj{\chi}(k') e^{i \pi \frac{k'^2}{m} } \sum_{k=1}^{4m} e^{\pi k i \frac{\xi-k'}{m}} (1+e^{-i \pi (m+k)}).
\end{align}
 After carrying out geometric sums over 1 to $4m$, we obtain
\begin{align}
J_2&=\frac{1-e^{4 i \pi \xi}}{2 \conj{C}} \sum_{k'=1}^{m} \conj{\chi}(k') e^{i \pi \frac{k'^2}{m} } \left(\frac{1}{e^{i \pi \frac{k'-\xi}{m}}-1} - \frac{e^{-i \pi m}}{e^{i \pi \frac{k'-\xi}{m}}+1} \right) \nonumber \\
&= (1-e^{4 i \pi \xi}) \frac{C}{2m} \sum_{k'=1}^{2m} \conj{\chi}(k') e^{i \pi \frac{k'^2}{m} } \frac{1}{e^{i \pi \frac{k'-\xi}{m}}-1},
\end{align}
where we have used the fact that $C \conj{C} = m$.
Thus overall we obtain
\begin{align}
\sum_{n=1}^{2m} J(\xi,n) \chi(n) e^{-i \pi \frac{n^2}{m}} &= \frac{e^{ \frac{i \pi}{m} \xi^2}}{1-e^{-2i \pi \xi}} \sum_{n=1}^{m} \chi(n) e^{-2i \pi \frac{\xi n}{m}} + \frac{C}{2m} \sum_{k'=1}^{2m} \conj{\chi}(k') e^{i \pi \frac{k'^2}{m} } \frac{1}{e^{i \pi \frac{k'-\xi}{m}}-1},
\end{align}
which can be written as
\begin{align}
\int_{-1 \swarrow 0} e^{ \frac{i \pi}{m} x^2-\frac{2i \pi}{m}\xi x} W(x) dx &= \frac{2 \pi i}{1-e^{-2i \pi \xi}} \sum_{n=1}^{m} \chi(n) e^{-2i \pi \frac{\xi n}{m}} - \frac{C}{2m} e^{-\frac{i \pi}{m} \xi^2} \conj{W(\conj{\xi})}.
\label{eq:w1}
\end{align}
Note that integral along ${-1 \swarrow 0}$  can be moved to contour $0 \swarrow 1$ on account of $\chi(0)=0$. This completes the proof of identity \autoref{e1}. 

From here on, we follow Siegel \cite{siegel1943contributions} to derive the main result in \autoref{eq:dirichlet}. We multiply both sides of \autoref{eq:w1} by $\xi^{s-1}$ and integrate over $\xi$ from 0 to $\infty e^{-\frac{\pi}{4} i}$. This leads to LHS being 
\begin{align}
\int_{0}^{\infty e^{-\frac{\pi}{4} i} } \xi^{s-1} \left( \int_{-1 \swarrow 0} e^{ \frac{i \pi}{m} x^2-\frac{2i \pi}{m}\xi x} W(x) dx \right)  d\xi & = \left( \frac{m}{2 \pi}\right)^s e^{-\frac{i \pi}{2}s} \Gamma(s) \lambda(s).
\end{align}
The first term on the RHS of \autoref{eq:w1} after being multiplied by $\xi^{s-1}$ and integrated over $\xi$ from 0 to $\infty e^{-\frac{\pi}{4} i}$ leads to
\begin{align}
\int_{0}^{\infty e^{-\frac{\pi}{4} i} } \frac{\xi^{s-1}}{1-e^{-2i \pi \xi}} \sum_{n=1}^{m} \chi(n) e^{-2i \pi \frac{\xi n}{m}} d\xi &= \int_{0}^{\infty e^{-\frac{\pi}{4} i} } \xi^{s-1} \sum_{n=1}^{\infty} \chi(n) e^{-2i \pi \frac{n}{m} \xi} d\xi \nonumber \\
&= \Gamma(s) \sum_{n=1}^{\infty} \chi(n) \left( \frac{2 i \pi n}{m} \right)^{-s} \nonumber \\
&= \left( \frac{m}{2 \pi}\right)^s e^{-\frac{i \pi}{2}s} \Gamma(s) L(s).
\end{align}
So \autoref{eq:w1} after being multiplied by $\xi^{s-1}$ and integrated over $\xi$ from 0 to $\infty e^{-\frac{\pi}{4} i}$ leads to
\begin{align}
 \left( \frac{m}{2 \pi}\right)^s e^{-\frac{i \pi}{2}s} \Gamma(s) \left(L(s)-\lambda(s)\right) & = C \int_{0}^{\infty e^{-\frac{\pi}{4} i} } e^{ -\frac{i \pi}{m} \xi^2} \xi^{s-1} \conj{W(\conj{\xi})} d\xi.
\end{align}
Using $W(-x)=-\chi(-1) W(x)$ which can be easily derived using $\chi(-1)\chi(n)=\chi(-n)$, we obtain
\begin{align}
 \left( \frac{m}{2 \pi}\right)^s e^{-\frac{i \pi}{2}s} \Gamma(s) \left(L(s)-\lambda(s)\right) & = \frac{C \conj{\lambda(1-\conj{s})}}{1+\chi(-1) e^{i \pi s}}.
\end{align}
With $a$ defined in \autoref{eq:a} and using the relation \cite{siegel1943contributions} 
\begin{equation}
2^s e^{\frac{\pi i}{2} s} \Gamma\left( \frac{s+a}{2}\right) \frac{1}{\Gamma(s)}= \pi^{-\frac{1}{2}} i^a \left( 1+\chi(-1) e^{i \pi s}  \right) \Gamma\left( \frac{1-s+a}{2}\right),
\end{equation}
we can finally obtain the desired result in \autoref{eq:dirichlet}.

\section{Double Exponential Method}
\label{sec:de}
The DE method \cite{takahasi1974double, mori2001double, trefethen2014exponentially}, also known anas $\text{tanh-sinh}$ method, is a method to compute certain numerical integrals where the integrands decay as double exponentially toward the boundary points. It uses the well known Simpson's quadrature method to calculate integrals to an astonishing accuracy of up to thousands of decimal places. We will adapt the method to work with oscillatory integrals. The integrand decays fast if the integral contour passes through a saddle point along the steepest direction. We will make a suitable variable change so that the integrand decays double exponentially away from the saddle point. We will introduce a further adjustment into this method in that nearby singularities of the transformed integrand will be included in estimating the integral.  

The ideas that we use are built upon the work of Turing \cite{turing1953some}, and Galway \cite{galway2001computing}. More recently, one can refer to \cite{al2021computation} which also discuss the idea of removing the contribution of singularities around the contour of integration.

Let us assume that we want to numerically calculate the integral of a function $f(z)$ along a contour $R$ that passes through the saddle point of $f(z)$ at $z_0$ along the steepest descent direction:
\begin{equation}
I=\int_{R}f\left( z\right) dz.
\label{eq:def_I}
\end{equation}
We would like to know the error $\Delta I_{h}=I-I_h$ that is made when the integral $I$ is
approximated by
\begin{equation*}
I_{h}=h\sum_{m=-\infty }^{\infty }f\left( z_{0}+mh\right) ,
\end{equation*}%
where $h$ represents a small complex interval chosen to discretise the contour $R$.  To get the error estimate it is helpful to define another integral, $I_C$,  along a closed contour $C=C_{U}+C_{L}$  with sides $C_{U}$ and $C_{L}$ parallel to $R$:
\begin{equation}
I_C = \int_{C}f\left( z\right) \phi _{L}\left( z\right) dz,
\label{eq:def_IC}
\end{equation}
where
\begin{equation}
\phi _{L}\left( z\right)=-\frac{1}{1-e^{2i \pi \frac{z-z_{0}}{h}}}.
\label{eq:phiL1}
\end{equation}
The integral $I_C$ is related to $I$ and $I_h$, and it will allow us to obtain an estimate of $\Delta I_{h}$ as will be shown below. We assume that the contribution of the side contours, $C_{\pm \infty}$ at $R \to \pm \infty$ is zero. This condition will be easy to meet in the problems that we discusses in this paper.
\begin{figure}[ht]
\centering
\includegraphics[width=0.9\textwidth]{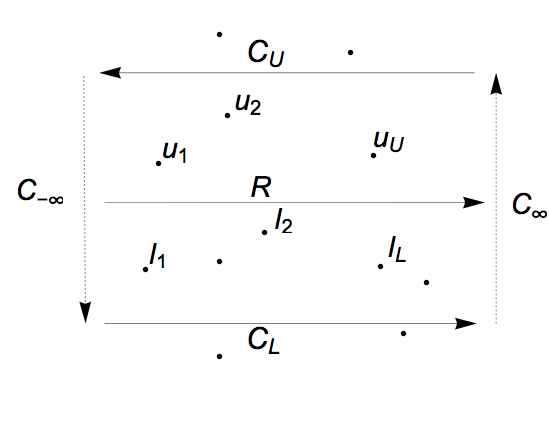}
\caption{Integral is to evaluated along the contour $R$. We consider the integral along $C_L \cup C_U \cup C_{-\infty} \cup C_{\infty}$. Contribution along $C_{\pm \infty}$ is assumed to be zero. It is assumed $f(z)$ has singularities between $R$ and $C_L$ ($C_U$) which are denoted by $l_i$ ($u_i$). One of these singularities lies near $C_L$ ($C_U$)}
\label{fig:contour}
\end{figure}
Note that $\phi_L(z)$ above is chosen such that it has poles at $z=z_0 + m h$ with residue $h$ where $m \in N$. The function $f(z)$ is assumed to have its own poles within $C$ which we collectively denote  by $\left\{ l_{0,}l_{1},\cdots l_{L},u_{0},u_{1},\cdots u_{U}\right\} $ where $l_{i}$ with $i\in \{0,1,..L\}$ denote the poles of $f(z)$ between $R$ and $C_{L}$ and $u_{i}$ with $i\in \{0,1,..U\}$ denote the poles of $f(z)$ between $R$ and $C_{U}$.
Application of the Residue theorem to $I_C$ leads to
\begin{align*}
\frac{I_{C}}{2i \pi}&=&\sum_{m=-\infty}^{\infty} \text{Res}[f(z)\phi _{L}, z_0+m h] +\sum_{i=0}^{U}\text{Res}[f(z)\phi _{L}(z),u_{i}] + \sum_{i=0}^{L}\text{Res}[f(z)\phi_{L}(z),l_{i}]\cdot 
\end{align*}
where Res$[g(z),a]$ denotes the residue of $g(z)$ at point $z=a$ and we assuming that $f(z)$ is not singular at $R$. Thus we have 
\begin{align*}
\frac{I_{C}}{2i \pi}&=&\frac{h}{2i \pi}\sum_{m=-\infty}^{\infty}f(z_0+m h) +\sum_{i=0}^{U}\text{Res}[f(z)\phi _{L}(z),u_{i}] + \sum_{i=0}^{L}\text{Res}[f(z)\phi_{L}(z),l_{i}]\cdot,
\end{align*}
which leads to 
\begin{align*}
\frac{I_{h}}{2 i \pi}&=\frac{I_C}{2 i \pi}-\sum_{i=0}^{U}%
\text{Res}[f(z)\phi _{L}(z),u_{i}]-\sum_{i=0}^{L}\text{Res}[f(z)\phi
_{L}(z),l_{i}]\cdot 
\end{align*}%
So the error $\Delta I_{h}=I-I_h$ is given by
\begin{align}
\frac{\Delta I_{h}}{2 i \pi} &=\int_{R}f\left( z\right) dz-\int_{C}f\left( z\right) \phi _{L}\left(
z\right) dz+\sum_{i=0}^{U}\text{Res}[f(z)\phi _{L}(z),u_{i}]+\sum_{i=0}^{L}\text{Res}[f(z)\phi _{L}(z),l_{i}],
\label{eqn:deltah}
\end{align}%
where we have used the definition of $I$ and $I_C$ from \autoref{eq:def_I} and \autoref{eq:def_IC}.
Now looking at closed contour formed from $R$, $C_{U}$ and irrelevant side contours at  $\pm \infty$, as shown in \autoref{fig:contour}, we can write
\begin{equation}
\int_{R}f\left( z\right) dz + \int_{C_{U}}f\left( z\right) dz =2i \pi~\sum_{i=0}^{U}\text{Res}[f(z),u_{i}],
\label{eqn:deltah1}
\end{equation}%
where Res$[f(z),u_{i}]$ refers to the residues of $f\left( z\right) $ at $u_{i}$ which lie between $R$ and $C_{U}$. 

Substituting the value of $\int_{R}f\left( z\right) dz$ from \autoref{eqn:deltah1} in \autoref{eqn:deltah} we obtain
\begin{align}
\frac{\Delta I_{h}}{2 i \pi}& =-\int_{C_{U}}f\left( z\right) dz-\int_{C_{U}+C_{L}}f\left(
z\right) \phi _{L}\left( z\right) dz +\sum_{i=0}^{U}\text{Res}[f(z)\phi _{L}(z),u_{i}]  \nonumber \\
&+\sum_{i=0}^{L}%
\text{Res}[f(z)\phi _{L}(z),l_{i}]+\sum_{i=0}^{U}\text{Res}[f(z),u_{i}] \nonumber \\
& =-\int_{C_{U}}f\left( z\right) \phi_{U}\left( z\right)
dz-\int_{C_{L}}f\left( z\right) \phi _{L}\left( z\right) dz \nonumber \\
& +\sum_{i=0}^{U}\text{Res}[f(z)\phi _{U}(z),u_{i}]+\sum_{i=0}^{L}%
\text{Res}[f(z)\phi _{L}(z),l_{i}],
\label{eq:I}
\end{align}%
where we have defined 
\begin{align}
\phi _{U}\left( z\right) &=\phi_{L}\left( z\right) +1 \nonumber \\
&=\frac{1}{(1-\exp \left(-2i \pi\left( z-z_{0}\right)/h\right) )}.
\end{align}
We can write \autoref{eq:I} expressions compactly as%
\begin{eqnarray}
I &=& I_{h}  + 2i \pi~\sum_{\{p_{i}\}}\text{Res}[f(z)\phi(z) ,p_{i}]  + \Delta I_{h}^{\prime},
\label{eq:int_error}
\end{eqnarray}%
where
\begin{equation*}
\Delta I_{h}^{\prime }=-2 i \pi \int_{C_{L}}f\left( z\right) \phi_{L}\left( z\right)
dz-2 i \pi \int_{C_{U}}f\left( z\right) \phi _{U}\left( z\right) dz.
\end{equation*}%
and $\left\{ p_{i}\right\} $ represent all singular points of $f(z)$ enclosed in $C$:  $\left\{ l_{0,}l_{1},\cdots l_{L},u_{0},u_{1},\cdots u_{U}\right\} $ and $\phi(z)$ defines a new function such that $\phi(z)=\phi_{L}(z)$ between $R$ and $C_L$ and  $\phi(z)=\phi_{U}(z)$ between $R$ and $C_U$.
It is clear from \autoref{eq:int_error} that a better representation of $I$ is given by
\begin{eqnarray}
I & \approx & I_{h}  + 2i \pi~\sum_{\{p_{i}\}}\text{Res}[f(z)\phi(z) ,p_{i}],
\label{eq:int_error_better}
\end{eqnarray}%
which is one of the central result of this paper.  The error in above approximation can be deduced as follows. We note that along $C_{L}$ the absolute value of $\exp \left( 2\pi
i\left( z-z_{0}\right) /h\right) $ is very large due to small $h$ and $%
\left( z-z_{0}\right) $ having a negative imaginary part and similarly along 
$C_{U}$ the absolute value of $\exp \left( -2i \pi\left( z-z_{0}\right)
/h\right) $ is very large due to $z$ having a positive imaginary part. This
leads to $\phi \left( z\right) $ being very small on $C.$ In addition, if the
path $C_{U}$ and $C_{L}$ are pushed out so that they pass near the other
singularities of $f(z)$ that have not been enclosed in $C$ then, we can approximate
these integrals using the saddle point method and thus obtain an approximation for $%
\Delta I_{h}^{\prime }$ . These saddle points develop due to an increase in $%
f\left( z\right) $ near the singularity and fast decay due to the
exponential factors in the denominators. Assuming singularity $z_{U}$ of $%
f\left( z\right) $ just outside $C_{U}$ and singularity $z_{L}$ of $f\left( z\right) 
$ just outside of $C_{L}$ we have
\begin{align}
\Delta I_{h}^{\prime }& \approx  -2 i \pi e^{- 2 i \pi \frac{z_{L}-z_{0}}{h}} \text{Res}[f(z)],z_{L}] + 2i \pi e^{2 i \pi \frac{z_{U}-z_{0}}{h}} \text{Res}[f(z),z_{U}] \label{error}
\end{align}

Note that we have used simple Simpson's rule for discretization. In principle, one can use the midpoint rule as well and define
\begin{equation*}
I_{h}=h\sum_{m=-\infty }^{\infty }f\left( z_{0}+ \left(m +\frac{1}{2}\right) h\right).
\end{equation*}%
In this case we will have
\begin{equation}
\phi _{L}\left( z\right)= \frac{1}{1+e^{2i \pi \frac{z-z_{0}}{h}}},
\label{eq:phiL2}
\end{equation}
which has poles at $z= z_0 + (m + 1/2) h$ with residue $h$ with $m \in N$. The rest of the treatment remains as above.  

Note that if there is a singularity very close to the contour of integration, then choosing simple Simpson's rule will result in a significant contribution from the node closest to the singularity. However, the correction term representing the removal of the singularity will be significant and effectively cancel the contribution from the node. The correction term from the singularity is significant because it is very close to the integration contour, and the $\phi_L$ becomes very large. To handle such cases where the singularities of the integrand fall close to the contour of integration, one may use the middle point rule. It works better because now the quadrature nodes can be put so that singularity lies midway between the two quadrature points, which ensures no node contributes overwhelmingly. Similarly, the contribution of the correction term does not blow up as the singularity approaches the contour of integration because now the correction term has $\phi _{L}$ is given by \autoref{eq:phiL2} which does not blow up in the limit as mentioned above. This difference in $\phi _{L}$ contrasts to $\phi _{L}$ factor as given by \autoref{eq:phiL1} which is valid for simple Simpson's rule.

\section{Algorithm}
\label{sec:algo}

\subsection{Riemann Zeta Function:}
\label{sec:algorithm}

We start with $I_0$ defined in \autoref{IN} and rewrite it as a series sum and a residual integral as follows. We define
\begin{eqnarray}
I_{N}\left( s\right)  
&=&\int_{N\swarrow (N+1)}\frac{e^{i\pi z^{2}}}{e^{i\pi z} -e^{-i\pi z}} z^{-s} dz, \label{IN}
\end{eqnarray}%
with integration performed along a straight line going from third quarter to
the first quarter and intersecting the x-axis between $N$ and $%
\left( N+1\right) $.  From the Residue theorem we have $%
I_{N}\left( s\right) =I_{\left( N-1\right) }\left( s\right) -N^{-s}$ which
leads to 
\begin{equation}
I_{0}\left( s\right) =\sum_{n=1}^{N}\frac{1}{n^{s}}+I_{N}\left( s\right) .
\label{eq:I0}
\end{equation}%
The reason to write $I_{0}\left( s\right)$ in terms of $I_{N}\left( s\right)$ is that that for a particular value
of $N$, that depends on $s$,  the function $I_{N}\left( s\right)$ develops a
saddle point and this allows for the application of the DE quadrature to numerically evaluate it.  In our case, we note that $\exp \left(
i \pi z^{2}\right) z^{-s}$ has saddle points at $z_{\pm}=\pm \sqrt{s/(2i \pi)}$.  Only the saddle point $z_{+}$ will be relevant for our case.  The direction of
steepest descent at $z_{+}$ is $\exp (i\pi /4).$ Thus we can transform $z$ along the steepest direction as
\begin{equation}
z\left( x\right) = r +\alpha \varepsilon \sinh\left(x\right)
\label{eq:transform}
\end{equation}
where $r = \sqrt{s/(2i \pi)}$ is the saddle point, $\varepsilon =$exp$\left(
i\pi /4\right) $ and $\alpha$ is a scaling factor that controls the distribution of singularities in the $x$ space.  A factor of $\sinh\left( x\right) $ is chosen to ensure that the integrand shows a double exponential decay along $x.$  We can further approximate the contour of integration $z\left( x\right)$ with 
$z(x) =N+1/2+\alpha \varepsilon \sinh(x)$ with $N=\lfloor 
\sqrt{t/(2\pi )}\rfloor $ where $\lfloor x\rfloor $ denotes the integer part
of $x.$ But a better choice of $N$ is given by
\begin{equation}
N =  \text{floor} \left(\Re{(r)} - \Im{(r)} \right).
\end{equation}
With the transform \autoref{eq:transform},  the integral \autoref{IN} can be written as
\begin{eqnarray}
I_{N}\left( s\right)  &=&\alpha \varepsilon \int_{-\infty }^{+\infty }f\left(
x\right) dx  \label{INX} \\
&=& -\alpha \varepsilon \int_{-\infty }^{+\infty }\frac{e^{i\pi z(x)^2}z(x)^{-s}}{e^{i \pi z(x)} -e^{-i\pi z(x)}}  \cosh(x) dx.
\end{eqnarray}%
Let us consider the location of singularities in $x$ space.
Singularities are defined by $\sin \left( \pi z\left( x\right) \right) =0$ which
after labelling the location of singularities with $x_{N,k}$ leads to $z\left(
x_{N,k}\right) =N+k$.  The $x_{N,m}$ thus satisfies 
\begin{equation}
r  + \alpha \varepsilon \sinh(x_{N,k}) = N + k
\end{equation}
Note that $x_{N,k}$ does not depend on $N$ and the use of $N$ here just indicates that $x_{N,k}$ is relative to $N$.  Note that the contour does not need to pass through the saddle point exactly. As for large $t$ the saddle point lies very close to $N$, one can choose $r$ to be $N+1/2$ in expansion \autoref{eq:transform}. Assuming this, the nearest singularities to the x-axis are given by
\begin{equation}
x_{N,k}=\ln \left( \frac{(k-1/2)}{\alpha \varepsilon } + \sqrt{\left( \frac{%
(k-1/2)}{\alpha \varepsilon }\right) ^{2}+1}\right) +2 i \pi n,~~~n \in \mathbb{Z}
\label{xm}
\end{equation}%
In particular, the two nearest singularities correspond to $k=0$ and $k=1$.  To provide an insight on these singularities we plot them for $\alpha =1$
and $\alpha =1/4$ for + sign above. The other singularities are displaced
with respect to these by $\pm 2n i \pi.$
\begin{figure}[ht]
\begin{minipage}[b]{0.47\linewidth}
\centering
\includegraphics[width=\textwidth]{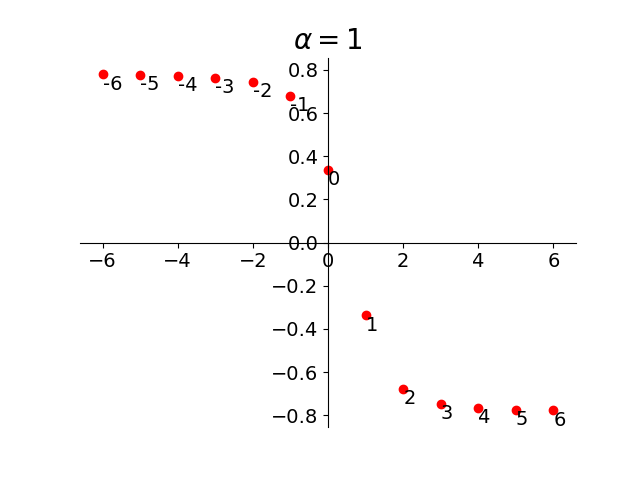}
\label{fig:figure1}
\end{minipage}
\hspace{0.5cm}
\begin{minipage}[b]{0.47\linewidth}
\centering
\includegraphics[width=\textwidth]{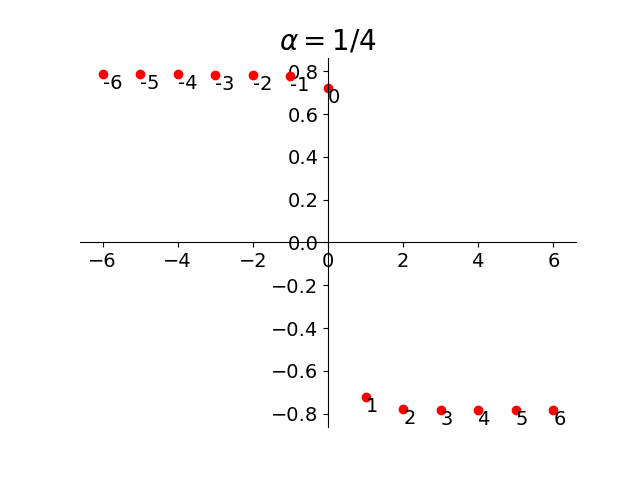}
\label{fig:figure2}
\end{minipage}
\caption{The distribution of singularities for $\alpha=1/4$ and $\alpha=1$ is shown. For smaller $\alpha$ singularities are pushed away from the $x$-axis to $y=\pi/4$. The singularities of the original integral at $N+k$ with $k>0$ are mapped to points with negative imaginary parts and those with $k\le 0$ mapped to the points with positive imaginary parts.}
\end{figure}
Choosing a smaller value of $\alpha$ pushes the
singularities away from the real axis. Note that asymptotically the location of singularities reaches $\pi/4$ for large $k$ values. 

\begin{figure}[ht]
\centering
\includegraphics[width=0.9\textwidth]{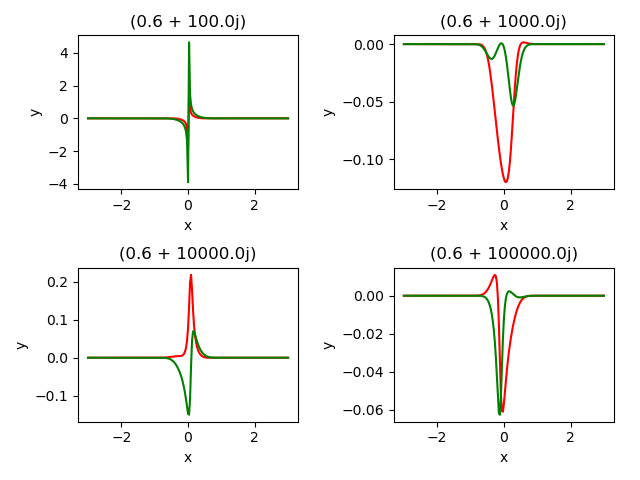}
\caption{Plot of real (red) and imaginary (green) values of the integrand for various $s$ values around the correct value of the saddle point. Note that integrands are not smooth in some cases as the contour may pass very close to a singularity due to the sine function in the denominator. However, this is not an issue for the application of the DE method.}
\label{fig:integrand_riemann}
\end{figure}

\begin{figure}[ht]
\centering
\includegraphics[width=0.9\textwidth]{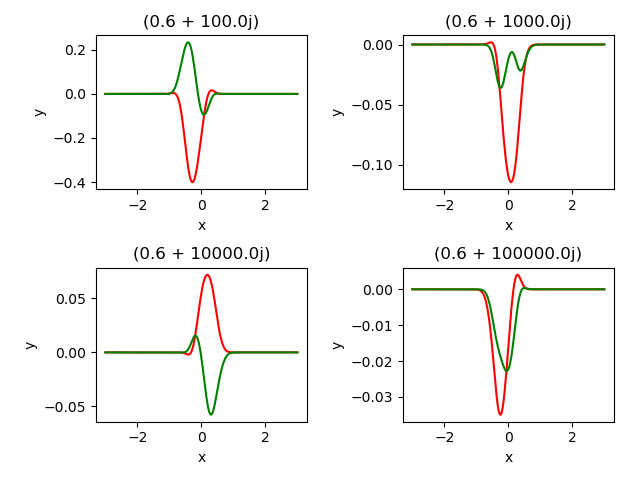}
\caption{Plot of real (red) and imaginary (green) values of the integrand for various $s$ values around the $r=n_0+1/2$.  This ensures the contour of integration is away from the singularities and this leads to relatively smooth plots.}
\label{fig:integrand_riemann_middle}
\end{figure}

A plot of integrand is shown in \autoref{fig:integrand_riemann} when the expansion is made around $r$ and in \autoref{fig:integrand_riemann_middle} when $r$ is chosen as $n_0+1/2$.

Applying the modified double exponential formula to this function $f\left( x\right) $ and
removing the contribution of all singularities within $C$ we obtain using \autoref{eq:int_error}
\begin{eqnarray}
I_{N}\left( s\right)  &=&\alpha \varepsilon \int_{-\infty }^{+\infty }\frac{e^{i\pi z(x)^2}z(x)^{-s}}{e^{i \pi z(x)} -e^{-i\pi z(x)}}  \cosh(x) \nonumber \\
&=&\alpha \varepsilon h\sum_{n=-\infty }^{\infty }\frac{e^{i\pi z(x_n)^2}z(x_n)^{-s}}{e^{i \pi z(x_n)} -e^{-i\pi z(x_n)}}  \cosh(x_n) + \sum_{x_{N,k}\in C}\frac{\phi \left( x_{N,k} \right) }{\left(
N+k\right) ^{s}}
\label{eq:riemann_sing}
\end{eqnarray}%
with approximate truncation error estimated from \autoref{error}.  In principle we can remove all singularities and write  \autoref{eq:I0} as:
\begin{eqnarray*}
I_{0}\left( s\right) 
&=&-\alpha \varepsilon h\sum_{n=-\infty }^{\infty }\frac{e^{i\pi z(x_n)^2}z(x_n)^{-s}}{e^{i \pi z(x_n)} -e^{-i\pi z(x_n)}}  \cosh(x_n) + \sum_{k=1}^{N}\frac{1}{k^{s}} \\
 &&- \sum_{k=1}^{\infty}\frac{\phi \left( x_{N,k-N} \right) }{k^{s}}
\end{eqnarray*}%
Note that $\phi=\phi_{L}$ for $k>N$ as the imaginary part of $x_{N,k}$ with $k \ge 0$ is negative and similarly $\phi=\phi_{U}$ for $m \le N$ and thus $\phi$ factors are quite small and thus the last summation represents a very minor correction to the overall result. The last two terms on the RHS can be combined and we can write the result compactly as
\begin{eqnarray*}
I_{0}\left( s\right) 
&=&-\alpha \varepsilon h\sum_{n=-\infty }^{\infty }\frac{e^{i\pi z(x_n)^2}z(x_n)^{-s}}{e^{i \pi z(x_n)} -e^{-i\pi z(x_n)}}  \cosh(x_n) - \sum_{k=1}^{\infty}\frac{\phi_{L} \left( x_{N,k-N} \right) }{k^{s}}.
\end{eqnarray*}%

Note that around the saddle point, the integrand decays double exponentially for large  $|x_{n}|$ as $\exp{\left( -2\pi \alpha ^{2}\sinh^2 \left( x\right)\right)}$. Thus the summation can be truncated at $x=\pm q$ where $q$ is chosen to
ensure that the required accuracy of $e^{-A}$ has been achieved. 
An estimate of $q$ is given by $\asinh \left( \sqrt{A/(2 \pi \alpha^2) }\right)$. We can approximate $h$ by equating the required accuracy with the error due to the contour integration. This leads to $h=\pi^{2}/(2A)$.  To be on the conservative side, we chose $h=\pi^{2}/(4A)$.  One can obtain a better estimate of $q$ and $h$ by root-finding, but the approximations above generally work well. Additionally, as we always work with extra precision of at least 20 decimal points, that ensures chosen $q$ and $h$ suffice to produce quite accurate results within the given error bounds.

\subsection{Lerch Function:}
\label{subsec:lerch}
\begin{figure}[ht]
\centering
\includegraphics[width=0.5\textwidth]{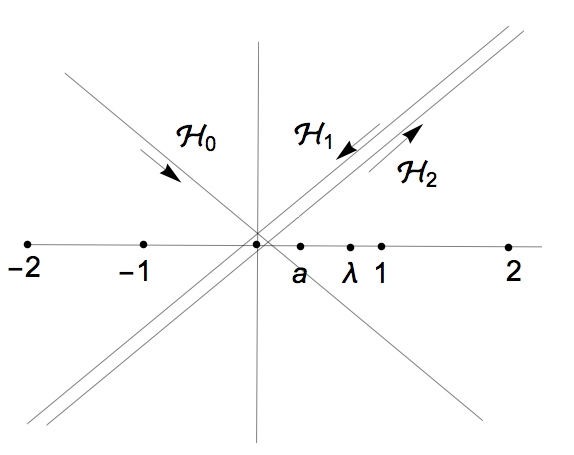}
\caption{ Hankel contour $\mathcal{H}$ is can be replaced with sum of two straight contours $\mathcal{H}_1$ and  $\mathcal{H}_2$.}
\label{fig:f2}
\end{figure}

The algorithm is based on \autoref{eq:lerch}. The contour $\mathcal{H}$ can be replaced with $\mathcal{H}_1$ and $\mathcal{H}_2$ as shown in \autoref{fig:f2}. On contour $\mathcal{H}_1$ it is important to assign angles to different points correctly. In particular,  if any point with $\theta > \pi/4$ with branch cut assumption at negative $x$-axis (the default assumption on most computer algebra systems) should be assigned an angle of $-(2\pi-\theta)$ if the branch cut is assumed at $\pi/4$. All three contours shown in \autoref{fig:f2} can be displaced along the $x$-axis to ensure RZ like expansions.

We denote the Lerch function as a combination of three integrations:
\begin{equation}
I = I_{\mathcal{H}_0} + I_{\mathcal{H}_1}+ I_{\mathcal{H}_2},
\end{equation}
where integrations on RHS are performed along contours $\mathcal{H}_0$,  $\mathcal{H}_1$ and $\mathcal{H}_2$. We consider all of them here integrations below.  We can write
\begin{align}
I_{\mathcal{H}_0}&=  e^{-i \pi a (1+a+2\lambda)}\int_{0 \searrow a} g_0(w) dw,
\label{eq:h0}
\end{align}
where
\begin{equation}
g_0(w) = \frac{e^{-i \pi w^2+ 2i \pi w (a+\lambda)}}{e^{\pi  i w- 2 i \pi a}-e^{-i \pi w}} w^{-s}.
\end{equation}
We consider the evaluation for the case where $t<0$. The $t>0$ case can be obtained from Lerch identity.   The integral has saddle points at the roots of $f'(w)=0$ where
\begin{align}
f(w)=w^2- 2 w (a+\lambda) + (s/(i \pi) \ln(w).
\end{align}
The saddle points correspond to 
\begin{equation}
w - (a+\lambda) + \frac{1}{w} \left( \frac{s}{2 i \pi}\right) = 0
\end{equation}
This leads to 
\begin{equation}
w = c \pm \sqrt{c^2- \frac{s}{2 \pi i}},
\end{equation}
where $c=(a+\lambda)/2$.
As there is a branch cut along the negative $x$ axis, we choose the saddle point on the positive $x$ and denote it by $w_0$:
\begin{equation}
w_0 \approx c + \sqrt{c^2- \frac{t}{2 \pi}}.
\end{equation}
Note that we have chosen $t<0$ and thus $w$ lies near the positive $x$ axis.  The $t>0$ case can be handled using the Lerch transformation \cite{LagariasWinnieLi+2012+1+48}:
\begin{equation}
L(s, \lambda, a) = \frac{\Gamma(s)}{(2\pi)^{1-s}}\left(e^{\frac{\pi i (1-s)}{2}} e^{-2\pi i \lambda a} L(1-s, 1-a, \lambda)+ e^{\frac{-\pi i (1-s)}{2}} e^{2 \pi i a (1-\lambda)} L(1-s, a, 1-\lambda) \right).
\end{equation}
The $t>0$ can also be handled by first calculating the conjugate of the Lerch function for appropriately transformed parameters and then take the conjugate of the result.

Note that for large $|t|$, the saddle point is approximately at $\sqrt{|t| / (2 \pi)}$.
We move the contour to pass through the saddle point. However it will encounter poles at $w = k + a$ with $k = 0,1,2,\cdots$.  The residue at $k$ is given by $e^{2 \pi i \lambda k} (k+a)^{-s}$.  We choose integer $n_0$ such that
\begin{equation}
n_0  = \text{floor}(w-a).
\end{equation}
Thus, we can write $I_{\mathcal{H}_0}$ as
\begin{equation}
I_{\mathcal{H}_0} = \sum_{k=0}^{n_0} \frac{e^{2 \pi i \lambda k}}{(k+a)^s} 
- e^{-i \pi a (1+a+2\lambda)}\int_{(n_0+a) \searrow (n_0+a+1)} g_0(w) dw.
\end{equation}
We will use the mapping
\begin{align}
w(x) &= w_0 + \frac{\alpha} {\varepsilon} \sinh(x).
\end{align}
The saddle points, $x_k$, near $w_0$ can be obtained by 
\begin{align}
w_0 + \frac{\alpha} {\varepsilon} \sinh(\tilde{x}_k) &= n_0 + a + k.
\end{align}
The change of variable leads to 
\begin{equation}
I_{\mathcal{H}_0} = \sum_{k=0}^{n_0} \frac{e^{2 \pi i \lambda k}}{(k+a)^s} 
+ \frac{\alpha e^{-i \pi a (1+a+2\lambda)}}{\varepsilon}\int_{-\infty}^{\infty}\ g_0(w(x)) \cosh(x) dx.
\end{equation}
The discretized version is given by
\begin{equation}
I_{\mathcal{H}_0} = \sum_{k=0}^{\infty} \frac{e^{2 \pi i \lambda k}}{(k+a)^s} \phi_U(\tilde{x}_{k-n_0})
+ \frac{h \alpha e^{-i \pi a (1+a+2\lambda)}}{\varepsilon}\sum_{n=-\infty}^{\infty} g_0(w(x_n)) \cosh(x_n).
\end{equation}
As with the RZ function, the $\phi_U(x_{n_0, k-n_0})$ factor is close to one before $n_0$ and close to zero after $n_0$. Unlike the $\phi_L$ factor for the RZ function, we have the $\phi_U$ factor in this case. The reason for this is the presence of $\varepsilon$ in the denominator rather than the numerator for Riemann zeta.

In the same manner we can write  $I_{\mathcal{H}_1}$ and  $I_{\mathcal{H}_2}$ as follows.  We note that the saddle points now correspond to $f'(w)=0$ where
\begin{align}
f(w)=w^2 + 2 w (a+\lambda) + (s-1)/(i \pi) \ln(w).
\end{align}
The saddle points are thus given by
\begin{equation}
w + (a+\lambda) + \frac{1}{w} \left( \frac{(s-1)}{2 \pi i}\right) = 0
\end{equation}
This leads to 
\begin{equation}
w_{\pm} = -c \pm \sqrt{c^2- \frac{(s-1)}{2 \pi i}},
\end{equation}
where $c=(a+\lambda)/2$. For large $|t|$, the saddle points fall at $\pm \sqrt{|t|/(2 \pi)}$. 

To calculate $I_{\mathcal{H}_1}$ we displace the contour $\mathcal{H}_1$ along the negative $x$ axis so that it passes through the saddle point, $w_{-}$, near $- \sqrt{|t|/(2 \pi)}$.   We choose positive integer $n_1$ such that
\begin{equation}
n_1  =- \text{floor}(w_{-}+\lambda).
\end{equation}

Taking appropriate poles into account this allows us to write 
\begin{equation}
\frac{I_{\mathcal{H}_1}}{q(s,\lambda)} = \int_{-\lambda \swarrow 0} g_1(w) dw,
\label{eq:h1}
\end{equation}
where 
\begin{equation}
g_1(w) = \frac{e^{i \pi w^2 +2 i \pi w(a+\lambda)}}{e^{i \pi w+2 \pi \lambda i}-e^{-i \pi w}} w^{s-1} ,
\end{equation}
and $q(s,\lambda)$ is pre-factor defined by
\begin{equation}
q(s,\lambda) =  i e^{i \pi \lambda (\lambda+1)+i \pi s/2} \frac{\Gamma(1-s)}{(2\pi)^{1-s}}.
\end{equation}
The integral above has poles at $(-k-\lambda)$ where we restrict integer values $k \ge 0$ so that brach cut at $\pi/4$ is not crossed.  The residue at pole $k$ is given by
\begin{align}
\text{Res}_k &=  \frac{e^{i \pi (-k-\lambda)^2 +2 i \pi (-k-\lambda)(a+\lambda)}}{2 \pi i e^{-i \pi (-k-\lambda)}} (-k-\lambda)^{s-1}  \\
&= -e^{-i \pi \lambda (\lambda+1)} e^{-2 \pi i a \lambda} \frac{e^{-i \pi s}}{2 \pi i} \frac{e^{-2 \pi i a k}}{(k+\lambda)^{(1-s)}} 
\end{align}

Taking contribution from all poles between the original and displaced contours and assuming $n_1>0$ we obtain,
\begin{equation}
I_{\mathcal{H}_1}= p(s, \lambda) \sum_{k=0}^{n_1-1} \frac{e^{-2 \pi i a k}}{(k+\lambda)^{(1-s)}} 
+ q(s, \lambda) \int_{-(n_1+\lambda) \nearrow -(n_1-1+\lambda)} g_1(w) dw,
\end{equation}
where
\begin{equation}
p(s,\lambda) =  e^{i \pi /2 (1-s)} e^{-2 \pi i a \lambda} \frac{\Gamma(1-s)}{(2\pi)^{1-s}}.
\end{equation}
For $n_1=0$ the contour passes between $-\lambda$ and 0 as in \autoref{eq:h1}.  The transformation used now is given by
\begin{align}
w &= w_{-} + \alpha \varepsilon \sinh(x),
\end{align}
and the saddle points near $w_{-} $ will be given by
\begin{align}
w_{-} + \alpha \varepsilon \sinh(\tilde{x}_k) &= -n_1 - k - \lambda
\end{align}
The discretized version is given by 
\begin{equation}
I_{\mathcal{H}_1}= p(s, \lambda) \sum_{k=0}^{\infty} \frac{e^{-2 \pi i a k}}{(k+\lambda)^{(1-s)}} \phi_L(\tilde{x}_{k-n_1})
+ hq(s, \lambda) \sum_{k=-\infty}^{\infty} g_1(w(x_k)) \cosh(x_n),
\end{equation}
To calculate $I_{\mathcal{H}_2}$ we displace the contour $\mathcal{H}_2$ along the positive $x$ direction so that the displaced contour passes close to the saddle point at $w_{+}$ which is approximately given by $+ \sqrt{|t|/(2 \pi)}$.  We choose positive integer $n_2$ such that
\begin{equation}
n_2  = \text{floor}(w_{+}+\lambda).
\end{equation}
Thus we have
\begin{equation}
\frac{I_{\mathcal{H}_2}}{q(s,\lambda)} = \int_{0 \nearrow (1-\lambda)} g_2(w)  dw,
\end{equation}
where
\begin{equation}
g_2(w) = \frac{e^{i \pi w^2 +2 i \pi w(a+\lambda)}}{e^{i \pi w+2 \pi \lambda i}-e^{-i \pi w}} w^{s-1}.
\end{equation}
The integrand above has poles at $(k-\lambda)$ where $k \in (1,2,\cdots) $ on the right side of the branch cut at $\pi/4$.  The residue at pole corresponding to $k$ is given by
\begin{align}
\text{Res}_k &=  \frac{e^{i \pi (k-\lambda)^2 +2 i \pi (k-\lambda)(a+\lambda)}}{2 \pi i e^{-i \pi (k-\lambda)}} (k-\lambda)^{s-1}  \\
&= e^{-i \pi \lambda (\lambda+1)} e^{-2 \pi i a \lambda} \frac{1}{2 \pi i} \frac{e^{2 \pi i a k}}{(k-\lambda)^{(1-s)}} 
\end{align}

Taking contribution from all poles between the original and displaced contours  and assuming $n_2>0$, we obtain
\begin{equation}
I_{\mathcal{H}_2}= r(s, \lambda) \sum_{k=1}^{n_2} \frac{e^{+2 \pi i a k}}{(k-\lambda)^{(1-s)}} 
- q(s, \lambda) \int_{(n_2-\lambda) \nearrow (n_2+1-\lambda)} g_2(w) dw,
\end{equation}
where
\begin{equation}
r(s,\lambda) =  e^{-i \pi /2 (1-s)} e^{-2 \pi i a \lambda} \frac{\Gamma(1-s)}{(2\pi)^{1-s}}.
\end{equation}
For $n_2=0$ the contour is confined between 0 and $1-\lambda$ only.  The transformation used now is given by
\begin{align}
w &= w_{+} + \alpha \varepsilon \sinh(x),
\end{align}
and the saddle points near $w_{+} $ will be given by
\begin{align}
w_{+} + \alpha \varepsilon \sinh(\tilde{x}_k) &= n_2 + k - \lambda.
\end{align}
The discretized version is 
\begin{eqnarray}
I_{\mathcal{H}_2}&=& r(s, \lambda) \sum_{k=1}^{n_2} \frac{e^{-2 \pi i a k}}{(k-\lambda)^{(1-s)}} \phi_L(\tilde{x}_{k-n_2}) \\
&&+ h~\alpha~\varepsilon~q(s, \lambda) \sum_{n=-\infty}^{\infty} g_2(w(nh)) \cosh(nh).
\end{eqnarray}

\subsection{Dirichlet $L$-Function:}
\label{subsec:dirichlet}
Our main implementation is based on \autoref{eq:dirichlet}.  We obtain the Hurwitz zeta function as a special case of the Lerch function and use
\autoref{eq:dirichlet_hurwitz} to obtain an alternative implementation of the function. Note that while the first implementation has a complexity of $\sqrt{q t}$, the second one has a complexity of $q \sqrt{t}$.
We write \autoref{eq:lambda} as 
\begin{align}
\lambda(x)=  \sum_{k=1}^{N}\frac{\chi(k)}{k^s} +  \frac{1}{2 \pi i} \int_{N \swarrow N+1} g(x) dx,
\label{eq:lambda1}
\end{align}
where 
\begin{equation}
g(x) = e^{\frac{i \pi}{m} x^2} x^{-s} W(x).
\end{equation}
The saddle point of the integrand lies at 
\begin{equation}
r = \sqrt{\frac{qs}{2 \pi i}}.
\end{equation}
Contour passing through the saddle point at an angle of $\pi/4$ will intersect between $N$ and $N+1$ where $N=\text{floor}(\Re(r)-\Im(r))$.
We define
\begin{align}
w &= r + \alpha \varepsilon \sinh(x),
\end{align}
and the saddle points near $w$ will be given by
\begin{align}
w + \alpha \varepsilon \sinh(\tilde{x}_k) &= N + k.
\end{align}
The discretized version is 

\begin{eqnarray}
\lambda(x) &=& \sum_{k=1}^{\infty}\frac{\chi(k)}{k^s} \phi_L(\tilde{x}_{k-N}) \\
&&+ h~\alpha~\varepsilon \sum_{n=\infty}^{\infty} g(w(nh)) \cosh(nh).
\end{eqnarray}

\section{Numeraical Experiments}
\label{sec:tests}

\subsection{Riemann Zeta Function:}
We test our algorithm for the RZ function against the one given by
Reyna \cite{de2011high, dereyna2022high} and implemented in MPMATH.  Both algorithms use the RS expansion for large $t$  and thus show an expected scaling as $\sqrt{t}$ for the evaluation of $\zeta (\sigma +it).$ For smaller values of $t$,  MPMATH uses other methods to evaluate the function while our algorithm still uses the RS expansion. Two important parameters of the algorithms are the $\alpha$ parameter and the number of singularities removed on either side of the integration contour. We choose $\alpha$ to be one and remove one singularity on either side of the contour by default. However when $|\Im{(s)}|$ is less than 100 we choose $\alpha$ to be 0.25. This lower $\alpha$ ensures that the branch cut pushes away from the contour. Next, we set the optimal cutoff parameter $q$ and step size $h$. We can either use a root-finding algorithm to find them or use the estimates derived in \autoref{sec:algo}. We use 20 extra decimal precision in internal calculations on top of desired precision. We carry out the tests in the critical strip, although the algorithm is valid even outside the critical strip. We choose values of $s=\sigma+i t$ with $\sigma \in \{0,0.1,\cdots,1\}$ and $t \in \{10^0,10^1,\cdots,10^8\}$. For each choice of $\sigma$ and $t$, we check the algorithm for different settings of decimal accuracy that is chosen from $\{10,20,100,200,400,500,1000\}$.  For all cases, $\zeta(s)$ matched perfectly well with MPMATH. This test shows that our algorithm is correct and that both implementations of the zeta function work well to arbitrary accuracy.

\begin{figure}[ht]
\begin{minipage}[b]{0.49\linewidth}
\centering
\includegraphics[width=\textwidth]{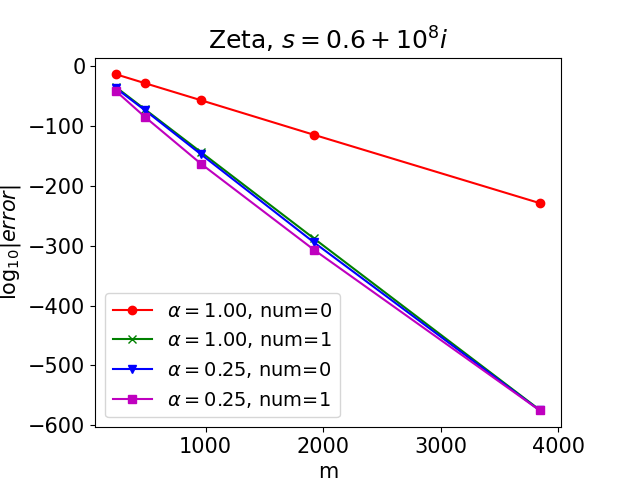}
\end{minipage}
\begin{minipage}[b]{0.49\linewidth}
\centering
\includegraphics[width=\textwidth]{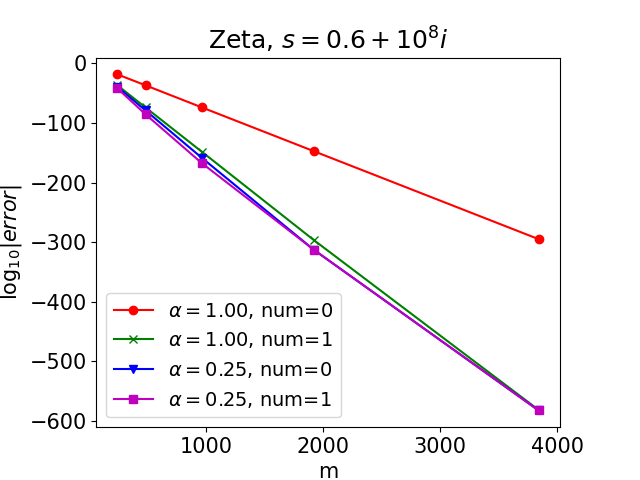}
\end{minipage}
\caption{Log error versus number of points used to discretise the residual integral ($m$). The figure on the left uses the saddle point as the middle point of discretisation grid while on the right uses the $(n_0+0.5)$ as the expansion point where $n_0$ is the nearest integer to the point where the integration contour intersects the x-axis. The parameter $\alpha$ is as in \autoref{eq:transform} and a smaller value of it pushes the singularities of integrand away from the x-axis. The 'num' parameter denotes the number of singularities removed from either side of the integration contour as per the second term on the RHS of \autoref{eq:riemann_sing}.}
\label{fig:riemann_h}
\end{figure}

Next, we show that the accuracy achieved for a given $q$ scales as $10^{-1/h}$ or equivalently as $10^{-m}$. In other words, we show that each halving of $h$ leads to a doubling of the accurate digits obtained. We take $s=0.6+10^8 i$  and set $q$ and $h$ parameters to obtain the RZ function with an accuracy of more than $10^{-600}$.  We take this as our reference value of the RZ function. In the subsequent experiment, we keep the $q$ parameter fixed to this value that results in an accuracy of at least $10^-{600}$. We then change $h \in \{1/2^1,1/2^2,\cdots,1/2^{5}\}$ and for each setting of $h$ we compute the zeta function. The error, $\epsilon$, is computed by taking the absolute value of the difference of $\zeta(s)$ for a given $h$ and the reference value. The results are shown in \autoref{fig:riemann_h} and they demonstrate a linear behaviour of $\ln(\epsilon)$ vs $m$. This is expected as for an appropriately chosen $q$, the total error in the zeta function computation is controlled by the size of $h$ and this error by \autoref{sec:algo} scales as $\ln(\epsilon)$ vs $m$. It is also clear from the test that only 4000 terms are enough to discretize the residual integral and obtain an accuracy of better than 500 digits. We note that this number of terms is quite small compared to terms in series on the RHS of \autoref{eq:I0} and thus contribute minimally to the overall computation cost of the algorithm. 

We now compare the timing of our implementation vs MPMATH. We note that there is a considerable scope of optimization. In particular, as the choice of $q$ and $h$ mainly depend on the accuracy desired, and not on the $t$, the calculation of the residual integral can be optimized by building a table of quadrature nodes and corresponding weights, and this table can then be used for all zeta function calculation corresponding to different $t$ values. We have not used these optimizations in the test. We choose $s=0.3 + 10^7 i$ for testing, set the accuracy to 400 decimal places, and compare the timing. Our implementation took about 2.3 seconds for it. This timing compares to 18.8 seconds taken by the MPMATH implementation. At 800 decimal accuracy, our timing was 12 seconds vs 146 seconds for the MPMATH. As we go for even higher $t$ and higher accuracy, MPMATH becomes significantly slower than our implementation.

\subsection{Lerch Function:}
As with the RZ function, when $|\Im(s)| \ge 100$, we work with $\alpha=1$ and remove the nearest singularity on either side of the integration contour. However, when $|\Im(s)|<100$, we work with $\alpha = 0.25$ and still remove one singularity on either side of the contour if possible. Note that in some cases, especially when $|\Im(s)|<100$, the integration contour can be close to the origin, and in such cases, there may be no singularity on the side closer to the origin. Other Lerch parameters are chosen as $a=0.3$ and $\lambda=0.7$.  

In the first experiment, we compare our implementation against MPMATH. We set $\Re(s) \in \{0, 0.1, 0.2, \cdots ,1 \}$ and $\Im(s) \in \{-1, -10, -100 \}$ and set our precision at 100 decimal points and use 150 decimal point precision for MPMATH which is enough for MPMATH to produce reliable values of the Lerch function for the chosen $s$ values here. A comparison of the results shows an excellent match with accuracy exceeding the preset limit of 100 in all cases. This test thus ensures that our implementation is correct.

In the second experiment, we check the timing. We set our precision at 50 decimal points and do the internal calculation with 20 decimals of extra precision for all $s$ cases considered below.   The precision set for our implementation is fixed at the value mentioned above and does not change. For the MPMATH function (LerchPhi), we set precision at different and typically higher values. High precision is needed because the MPMATH function produces unintuitive results at low precisions, especially for $t>100$. We consider our result the reference one and increase the MPMATH precision in steps of 10 until it shows the same value as our implementation with better than 50 decimal accuracy.   We first choose $s=0.6-100 i$. For this value of $s$, we need to set the MPMATH precision at 90, and it takes MPMATH 0.7 seconds vs 0.2 seconds for our implementation. Next we set $s=0.6-1000 i$.  For this value of $s$,  LerchPhi requires precision of 650 decimal points. It takes LerchiPhi 27 seconds while it takes 0.2 seconds for our implementation. Finally, we choose $s=0.6-10^8 i$ and obtain the following result for the Lerch function to 50 decimal accuracy with our implementation:
\begin{eqnarray}
(-1.642781971616430577623708444116579252538445671826 - \nonumber \\  4.4985920038844187868475845606656684449691018721039 i). \nonumber
\end{eqnarray}
It took about 0.78 seconds to obtain the above result. Unfortunately, both MPMATH and Mathematica implementations cannot handle cases like these.  

We check that the Lerch function degenerates to  Hurwitz zeta and RZ functions in appropriate limits. For this purpose, the RZ function as derived from the Lerch function was tested and found to be in excellent agreement with the implementation of the RZ function. as described in \autoref{sec:algorithm}. Next, we calculate the Dirichlet $L$-function through its connection to the Hurwitz zeta functions and compare the results with an alternative implementation of the function based on the Siegel formula as described in \autoref{subsec:dirichlet}. The results show an excellent agreement. The test thus shows that the Lerch zeta function implementation can reproduce the formulas for the RZ, the Hurwitz zeta function and the Dirichlet $L$-functions.

\begin{figure}[ht]
\begin{minipage}[b]{0.49\linewidth}
\centering
\includegraphics[width=\textwidth]{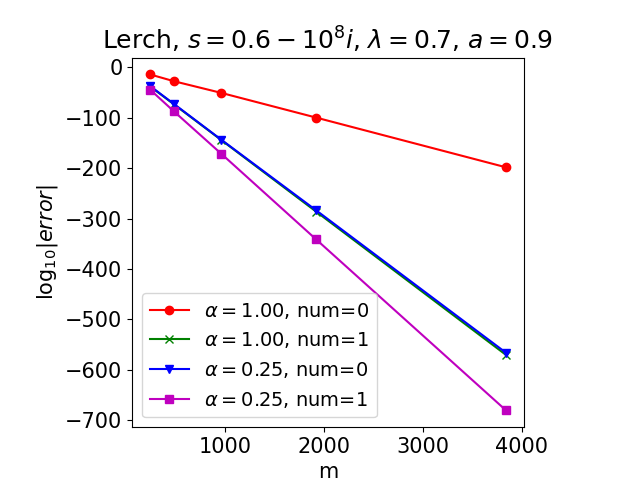}
\end{minipage}
\begin{minipage}[b]{0.49\linewidth}
\centering
\includegraphics[width=\textwidth]{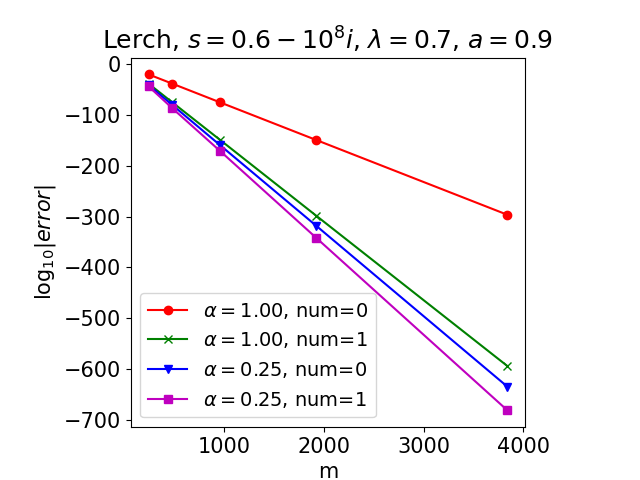}
\end{minipage}
\caption{Log error versus number of points used to discretise the residual integral ($m$). The figure on the left uses the saddle point as the central node of the integral discretisation while on the right we use $(n_0+0.5)$ as the central node where $n_0$ is the nearest integer to the point where the integration contour intersects the x-axis. The parameter $\alpha$ is as discussed in \autoref{subsec:lerch} and a smaller value of it pushes the singularities of integrand away from the x-axis. The 'num' parameter denotes the number of singularities removed from either side of the integration contour.}
\label{fig:lerch_h}
\end{figure}

Lastly, we check that the accuracy obtained as a function of $1/h$ or $m$ shows a linear scaling. The experiment is very similar to the one carried out for the RZ function. The results for $s=0.6-10^8 i$ are shown in \autoref{fig:lerch_h}. The results show that every halving of $h$ leads to a doubling of the precision. Further, the impact of removing the contribution of the singularities is clear. For $\alpha=1$, removing the contribution of nearby singularities, one on each side of the contour, leads to a significant improvement in the accuracy obtained. Further, setting $\alpha=0.25$ and not removing any contribution from the singularity leads to very accurate results too. In this case, we achieve higher accuracy because the singularities move away from the integration contour.

\subsection{Dirichlet L-Function:}
We work with Dirichlet character 
\begin{align*}
\chi_{8,2} &=[\chi(0), \chi(1), \cdots,\chi(7)] \\
 &= [0, 1, 0, -1, 0, -1, 0, 1]
\end{align*}
and set $s=0.6+ 10^5 i$. We compare our results against MPMATH for various accuracies.  At 50 decimal accuracy our results tie up exactly with MPMATH but while our implementation takes 1.1 second, the MPMATH implementation takes about 8.7 seconds.  At 200 decimal accuracy again the results tie up and time taken by our implementation is about 5.9 seconds vs 11.5 seconds for MPMATH.  We further considered higher accuracies of up to 500 and in all cases there was a perfect agreement between MPMATH and our implementations.  

We set $s=0.6+ 10^6 i$ and the desired accuracy to be 50. Our algorithm takes only 1.2 seconds, while MPMATH takes 71.7 seconds. As with all previous cases,  results agree better than the preset limit of 50 decimal points. We double the imaginary part of $s$ and take $s=0.6+2 \times 10^6$.  Now MPMATH takes about 144 seconds, while our implementation still takes about 1.2 seconds. Any further increase in $t$ or accuracy required leads to the MPMATH implementation being considerably slower than our implementation. The reason for the slowdown is that while our algorithm has a $\sqrt{t}$ complexity, the MPMATH algorithm has $t$ complexity.  

We now consider two cases where MPMATH is very slow and thus a comparison against it is not possible. For these cases, we use another implementation of Dirichlet $L$-function by writing it in terms of the Hurwitz zeta function. The Hurwitz zeta function itself is calculated as special case of the Lerch zeta function as discussed in \autoref{subsec:hurwitz}. We work with Dirichlet character 
\begin{align*}
\chi_{7,5} &=[\chi(0), \chi(1), \cdots,\chi(6)] \\
 &= [0, 1, w^2, -w, -w, w^2, 1],
 \end{align*}
where $w = \exp(i \pi/3)$ and set the accuracy goal to be 50 digits. We work with $s=0.6 + 10^{8} i$ and obtain the following result for the Dirichlet $L$-function from both implementations:
 \begin{eqnarray}
(0.34580337824253257378760299316255985284400906588262 - \nonumber \\  1.0760292785488344655945930583565334551150785880025 i). \nonumber
\end{eqnarray}
The Siegel Dirichlet $L$-function formula took about 1.2 seconds for the calculation above, while the Hurwitz zeta function-based approach took 12.6 seconds. We checked the results at 100 digits accuracy, and again the results matched, and while the first implementation took 1.9 seconds, it took 18.3 seconds with the second implementation.   At 500 digit accuracy, the first implementation took about 15.7 seconds vs 152.0 seconds for the second implementation. Similarly, working with $s=0.6 + 10^{10} i$, we obtained excellent agreement in results at various accuracies of 50, 100 and 500. Unfortunately, both MPMATH and Mathematica implementations cannot handle cases like these.

\begin{figure}[ht]
\begin{minipage}[b]{0.49\linewidth}
\centering
\includegraphics[width=\textwidth]{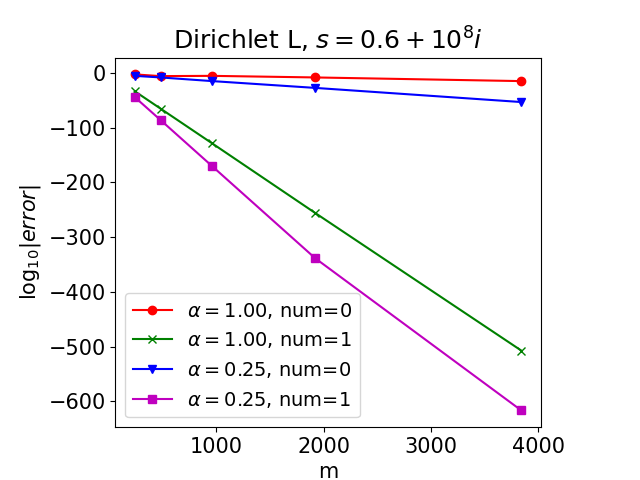}
\label{fig:dirichlet_h}
\end{minipage}
\begin{minipage}[b]{0.49\linewidth}
\centering
\includegraphics[width=\textwidth]{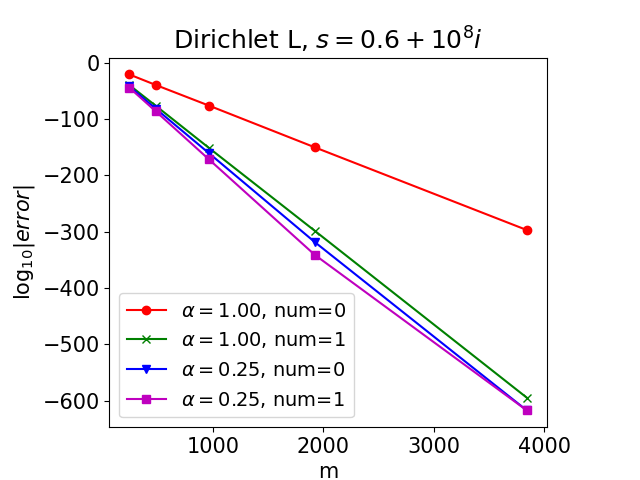}
\label{fig:dirichlet_h_n0}
\end{minipage}
\caption{Log error versus number of points used to discretise the residual integral ($m$). The figure on the left uses the saddle point as the central node of the integral discretisation while on the right we use $(n_0+0.5)$ as the central node where $n_0$ is the nearest integer to the point where the integration contour intersects the x-axis. The parameter $\alpha$ is as discussed in \autoref{subsec:lerch} and a smaller value of it pushes the singularities of integrand away from the x-axis. The 'num' parameter denotes the number of singularities removed from either side of the integration contour.}
\end{figure}

\section{Conclusion}
\label{sec:conclusion}
We have developed the MDE quadrature formulas to numerically calculate Riemann zeta, Lerch and Dirichlet $L$-functions to thousands of decimal point accuracy in the all-important critical region. Our method is based on the RS formulas and treats the remainder integrals exactly. The method retains all the good features of the underlying RS formulas such $\sqrt{t/(2\pi)}$ or $\sqrt{q t/(2\pi)}$ scaling. The extra computational effort needed to treat the residual integrals exactly is minimal. The method provides a very accurate estimate of the error. We have shown that the number of expected good decimal points doubles with every halving of the discretization size. Thus, the method provides a way to control the error automatically. Computing the value for a given $h$ and subsequently halving $h$ can indicate the number of digits one can trust. We have provided the numerical estimate of the quadrature parameters, and they produced the desired accuracy for all the cases tested in this paper. One can download the Python scripts showing the implementation of these functions from Github \cite{tyagi2022}. They allow one to reproduce almost all numerical results reported in this paper.

The author has successfully used the method to calculate many other functions including incomplete Gamma function for complex arguments \cite{temme1996special, gil2007numerical}, $K$-Bessel function of imaginary order \cite{booker2013bounds}, truncated theta sums \cite{hiary2011nearly} and many other functions.  We can calculate all of the above functions to thousands of digits of accuracy. We may report these results in future papers.

\bigskip 
\bibliographystyle{unsrt} 
\bibliography{riemann}

\end{document}